\documentclass[english]{smfartsMM}
\usepackage{smfthm}
\usepackage{MnSymbol}
\usepackage[unicode, pdftex]{hyperref}
\usepackage{epigraph}
\usepackage{tikz-cd}
\usepackage{parskip}
\usepackage[english]{babel}
\usepackage[utf8]{inputenc}
\usepackage[matrix, arrow,all,cmtip,color]{xy}
\usepackage{graphicx}

\def\Comment#1#2{#1}{}
\def\comment#1{}

\def\includegraphicss[#1]#2{ 
             \pdfximage width  #1\linewidth {#2} 
             \pdfrefximage\pdflastximage 
             }

\newtheorem{theorem}{Theorem}[section]

\def\Imm{\operatorname{Im}}
\def\cl{\operatorname{cl}}
\def\clB{\operatorname{cl_B}}

\def\Sets{\ensuremath{\mathrm{Sets}}}

\def\Topp{\ensuremath{\mathrm{Top}}}
\author{M.Gavrilovich and M.Rabivovich} 
\title
[The Quillen negation monoid]
{The Quillen negation monoid of a category, and Schreier graphs of its action on classes of morphisms}

\thanks{Misha Gavrilovich and Misha Rabinovich. The University of Haifa, Israel. \\ 
{\tt{mi\!\!\!ishap\!\!\!p@sd\!\!\!df.org}} and {\tt{stormspi\!\!\!irit1\!\!\!16@ma\!\!\!ail.ru}}  
The first author was supported by ISF grant 290/19.   
}

\email{mi\!\!\!ishap\!\!\!p@sd\!\!\!df.org}

\renewcommand\star{\filledstar}

\def\smallblackbox{{\Huge\centerdot}}
\def\bblackbox{\blacksquare}
\def\bblacksquareR{\resizebox{3.39pt}{!}{$\blacksquare$}}
\newcommand{\bblacksquare}{\mathpalette\raisebblacksquare\relax}
\newcommand{\raisebblacksquare}[2]{\raisebox{2\depth}{$#1\bblacksquareR$}}

\def\ptt{\{\bullet\}}
 
\def\lra{\longrightarrow} 
\def\xra{\xrightarrow} 
\def\rtt{\operatorname{\rightthreetimes}}
\def\inv{^{-1}}

\def\MtoL{ {\left\{\underset{\bblacksquare_A}{}{\swarrow} \overset{\raisebox{1pt}{$ \bullet_{U}$}}{}{\searrow} 
\underset{\bblacksquare_{X}}{}
{\swarrow}\overset{\raisebox{1pt}{$\bullet_{V}$}}{}{\searrow}\underset{\bblacksquare_B}{}\right\}
}\atop{\Downarrow\atop
{	 
\left\{\underset{\bblacksquare_A}{}{\swarrow} \overset{\raisebox{1pt}{$ \bullet_{U=X=V}$}}{}{\searrow} 
\underset{\bblacksquare_B}{}\right\}
}}
}

\def\mathcalDense{\protect{\hphantom{\bullet_c\searrow}{\!\smallblackbox \atop  \Downarrow }\atop {\{\raisebox{3pt}{$\bullet_o$}\searrow\raisebox{-3pt}{$\smallblackbox_c$}\} }}} 
\def\mathcalDenseSection{{\,\,\,\,\,\,{\!\smallblackbox_c \atop  \Downarrow }\atop {\{\raisebox{3pt}{$\bullet_o$}\searrow\raisebox{-3pt}{$\smallblackbox_c$}\} }}} 
\def\mathcalDensebrackets{\left(\mathcalDense\right)}
\def\nondenseimage{\mathcalDense} 
\def\nondenseimagebrackets{\mathcalDensebrackets}

\def\nondenseimageSectionbrackets{\left(\mathcalDenseSection\right)} 

\def\Otoc{\{\raisebox{2pt}{$\bullet$}\searrow\raisebox{-2pt}{$\smallblackbox$}\}}

\def\mathcalD{\protect{\{\raisebox{3pt}{$\bullet$}\searrow\raisebox{-3pt}{$\smallblackbox$}\}\atop{\Downarrow\atop\{\bullet\}}}}
\def\mathcalDbrackets{\left(\mathcalD\right)}

\def\mathcalE{{\{\star\leftrightarrow\star\}\atop{\Downarrow\atop\{\star\}}}}

\def\mathcalEbrackets{\left({\{\star\leftrightarrow\star\}\atop{\Downarrow\atop\{\star\}}}\right)}

\def\mathcalDE{\mathcalD,\mathcalEnobrackets}
\def\mathcalDEbrackets{\left(\mathcalDE\right)}

\def\mathcalDEmorphism{{\{\star\leftrightarrow\star\rightarrow\smallblackbox\}\atop{\Downarrow\atop\{\bullet\}}}}
\def\mathcalDEmorphism{{\protect{\{\raisebox{3pt}{$\star\leftrightarrow\star$}\searrow\raisebox{-3pt}{$\smallblackbox$}\}}\atop{\Downarrow\atop\{\bullet\}}}}
\def\mathcalDEmorphismbrackets{\left(\mathcalDEmorphism\right)}

\def\mathcalEDmorphism{{\{\bullet\rightarrow \star\leftrightarrow\star\}\atop{\Downarrow\atop\{\bullet\}}}}
\def\mathcalEDmorphismbrackets{\left(\mathcalDEmorphism\right)}

\def\mathcalDD{{\{\raisebox{3pt}{$\bullet$}\searrow\raisebox{-3pt}{$\smallblackbox$}\}\atop{\Downarrow\atop\resizebox{22pt}{!}{$\{\star\leftrightarrow\star\}$}}}}
\def\mathcalDDbrackets{\left(\mathcalDD\right)}

\def\mathcalDDSection{\protect{\{\protect\raisebox{3pt}{$\bullet$}\searrow\protect\raisebox{-3pt}{$\centerdot$}\}\atop{\Downarrow\atop\{\star\leftrightarrow\star\}}}}
\def\mathcalDDbracketsSection{\left(\mathcalDDSection\right)}

\def\lra{\longrightarrow}

\def\mathcalA{\left({{\{\bullet,\bullet\}}\atop{\overset{\Downarrow}{\{\bullet\}}}}\right)}
\def\mathcalB{\{\{\star\}\lra\{\star\leftrightarrow\star\}\}}
\def\mathcalB{\left(\displaystyle{\overset{\star}{\Downarrow}\,\,\,\,\,\,\,\,}\atop{\{\star\leftrightarrow\star\}}\right)}
\def\mathcalC{\left({\emptyset\atop{\Downarrow\atop\{\bullet\}}}\right)}
\def\mathcalCnobrackets{{\emptyset\atop{\Downarrow\atop\{\bullet\}}}}

\def\mathcalBnobrackets{\displaystyle{\overset{\star}{\Downarrow}\,\,\,\,\,\,\,\,}\atop{\{\star\leftrightarrow\star\}}}
\def\mathcalCnobrackets{{\emptyset\atop{\Downarrow\atop\{\bullet\}}}}

\def\mathcalA{\left({{\{\bullet,\bullet\}}\atop{\overset{\Downarrow}{\{\bullet\}}}}\right)}
\def\mathcalAnobrackets{{\{\bullet,\bullet\}\atop{\displaystyle \Downarrow \atop\ptt}}}

\def\mathcalB{\{\{\star\}\lra\{\star\leftrightarrow\star\}\}}
\def\mathcalB{\left(\displaystyle{\overset{\star}{\Downarrow}\,\,\,\,\,\,\,\,}\atop{\{\star\leftrightarrow\star\}}\right)}
\def\mathcalC{\left({\emptyset\atop{\Downarrow\atop\{\bullet\}}}\right)}
\def\mathcalCnobrackets{{\emptyset\atop{\Downarrow\atop\{\bullet\}}}}

\def\mathcalE{(\{\bullet\leftrightarrow\bullet\}\lra \{\bullet\})}
\def\mathcalE{\left({\{\star\leftrightarrow\star\}\atop{\Downarrow\atop\{\star\}}}\right)}
\def\mathcalEnobrackets{{\{\star\leftrightarrow\star\}\atop{\Downarrow\atop\{\star\}}}}
\def\mathcalF{{\!\!{\bullet\leftrightarrow\bullet_\downarrow\phantom{\smallblackbox\leftrightarrow\smallblackbox}}\atop{\phantom{\star\star}\smallblackbox\leftrightarrow\smallblackbox}}
	     \atop{\Downarrow\atop{\star\leftrightarrow\star\leftrightarrow\star}}}

\def\mathcalFwithIndices{{\!\!{\bullet_u\leftrightarrow\bullet_{\!\!v\atop\downarrow}\phantom{\smallblackbox\leftrightarrow\smallblackbox}}\atop{\phantom{\star_u\leftrightarrow\star_v}\smallblackbox_c\leftrightarrow\smallblackbox_d}}
	     \atop{\Downarrow\atop{\star_u\leftrightarrow\star_{v=c}\leftrightarrow\star_d}}}

\def\mathcalFbrackets{\left({\!\!{\bullet\leftrightarrow\bullet_\downarrow\phantom{\smallblackbox\leftrightarrow\smallblackbox}}\atop{\phantom{\star\star}\smallblackbox\leftrightarrow\smallblackbox}}
	     \atop{\Downarrow\atop{\star\leftrightarrow\star\leftrightarrow\star}}\right)}

\def\mathcalFbracketsSection{\left({\!\!{\bullet\leftrightarrow\bullet_\downarrow\,\,\,\,\,\,\,}\atop{\,\,\,\,\centerdot\leftrightarrow\centerdot}}
	     \atop{\Downarrow\atop{\star\leftrightarrow\star\leftrightarrow\star}}\right)}

\def\mathcalG{\left({\bullet\phantom{,\bullet}}\atop{{\!\!\!\Downarrow}\atop{\bullet,\bullet}}\right)}

\def\mathcalGnobrackets{{\bullet\phantom{,\bullet}}\atop{{\!\!\!\Downarrow}\atop{\bullet,\bullet}}}

\def\mathcalGnobracketsSection{{\bullet\,\,\,\,\,\,\,\,}\atop{{\!\!\!\Downarrow}\atop{\bullet,\bullet}}}

 \def\closedsubspace{{{\left\{\raisebox{2pt}{$\star_x\leftrightarrow\star_y\leftrightarrow\star_z$} \searrow \raisebox{-2pt}{$\smallblackbox_c$}\right\}}\atop{\Downarrow}}
 \atop{\{ \star_x=\star_y\leftrightarrow \star_z=\bblackbox_c\}}}

 \def\closedsubspace{{{\raisebox{2pt}{$\star_x\leftrightarrow\star_y\leftrightarrow\star_z$} \searrow \raisebox{-2pt}{$\smallblackbox_c$}\right\}}\atop{\Downarrow}}
 \atop{\displaystyle\{ \star_{x=y}\leftrightarrow \star_{z=c}\}}}

\def\MtoL{ { \underset{\bblacksquare_A}{}{\swarrow} \overset{\raisebox{1pt}{$ \bullet_{U}$}}{}{\searrow} 
\underset{\bblacksquare_{X}}{}
{\swarrow}\overset{\raisebox{1pt}{$\bullet_{V}$}}{}{\searrow}\underset{\bblacksquare_B}{} 
}\atop{\displaystyle\Downarrow\atop
{	 
 \underset{\bblacksquare_A}{}{\swarrow} \overset{\raisebox{1pt}{$ \bullet_{U=X=V}$}}{}{\searrow} 
\underset{\bblacksquare_B}{} 
}}
}

\def\mathcalDense{\protect{\hphantom{\bullet_c\searrow}{\!\smallblackbox \atop  \displaystyle\Downarrow }\atop { \raisebox{3pt}{$\bullet$}\searrow\raisebox{-3pt}{$\smallblackbox$}  }}} 
\def\mathcalDenseSection{{\,\,\,\,\,\,{\!\smallblackbox \atop  \displaystyle\Downarrow }\atop { \raisebox{3pt}{$\bullet$}\searrow\raisebox{-3pt}{$\smallblackbox$}  }}} 
\def\mathcalDensebrackets{\left(\mathcalDense\right)}
\def\nondenseimage{\mathcalDense} 
\def\nondenseimagebrackets{\mathcalDensebrackets}

\def\nondenseimageSectionbrackets{\left(\mathcalDenseSection\right)}

\def\OtoOc{{{{\displaystyle \bullet\phantom{\searrow\!\smallblackbox}} \atop  \displaystyle\Downarrow }\atop { \raisebox{3pt}{$\bullet$}\searrow\raisebox{-3pt}{$\smallblackbox$}  }}}

\def\Otoc{ \raisebox{2pt}{$\bullet$}\searrow\raisebox{-2pt}{$\smallblackbox$} }

\def\mathcalD{\protect{ \raisebox{3pt}{$\bullet$}\searrow\raisebox{-3pt}{$\smallblackbox$} \atop{\displaystyle\Downarrow\atop \displaystyle \bullet }}}
\def\mathcalDbrackets{\left(\mathcalD\right)}

\def\mathcalE{{ \star\leftrightarrow\star \atop{\displaystyle\Downarrow\atop \star }}}

\def\mathcalEbrackets{\left({ \star\leftrightarrow\star \atop{\displaystyle\Downarrow\atop \star }}\right)}

\def\mathcalDE{\mathcalD,\mathcalEnobrackets}
\def\mathcalDEbrackets{\left(\mathcalDE\right)}

\def\mathcalDEmorphism{{\protect{ \raisebox{3pt}{$\star\leftrightarrow\star$}\searrow\raisebox{-3pt}{$\smallblackbox$} }\atop{\displaystyle\Downarrow\atop \displaystyle \bullet }}}
\def\mathcalDEmorphismbrackets{\left(\mathcalDEmorphism\right)}

\def\mathcalEDmorphism{{\protect{ \raisebox{3pt}{$\bullet$}\searrow\raisebox{-3pt}{$\star\leftrightarrow\star$} }\atop{\displaystyle\Downarrow\atop \displaystyle \bullet }}}
\def\mathcalEDmorphismbrackets{\left(\mathcalEDmorphism\right)}

\def\mathcalDD{{ \raisebox{3pt}{$\bullet$}\searrow\raisebox{-3pt}{$\smallblackbox$} \atop{\displaystyle\Downarrow\atop\resizebox{22pt}{!}{$ \star\leftrightarrow\star $}}}}
\def\mathcalDDbrackets{\left(\mathcalDD\right)}

\def\mathcalDDSection{\protect{ \protect\raisebox{3pt}{$\bullet$}\searrow\protect\raisebox{-3pt}{$\smallblackbox$} \atop{\displaystyle\Downarrow\atop \star\leftrightarrow\star }}}

\def\mathcalDDbracketsSection{\left(\mathcalDDSection\right)}

\def\lra{\longrightarrow}

\def\mathcalA{\left({{ \bullet,\bullet }\atop{\overset{\displaystyle\Downarrow}{ \bullet }}}\right)}
\def\mathcalB{  \star \lra \star\leftrightarrow\star  }
\def\mathcalB{\left(\displaystyle{\overset{\star}{\displaystyle\Downarrow}\,\,\,\,\,\,\,\,}\atop{ \star\leftrightarrow\star }\right)}
\def\mathcalC{\left({\emptyset\atop{\displaystyle\Downarrow\atop \displaystyle \bullet }}\right)}
\def\mathcalCnobrackets{{\emptyset\atop{\displaystyle\Downarrow\atop \displaystyle \bullet }}}

\def\mathcalBnobrackets{\displaystyle{\overset{\star}{\displaystyle\Downarrow}\,\,\,\,\,\,\,\,}\atop{ \star\leftrightarrow\star }}
\def\mathcalCnobrackets{{\emptyset\atop{\displaystyle\Downarrow\atop \displaystyle \bullet }}}

\def\mathcalA{\left({{\displaystyle  \bullet\,\,\bullet }\atop{\overset{\displaystyle\Downarrow}{ \displaystyle \bullet }}}\right)}
\def\mathcalAnobrackets{{\displaystyle \bullet\,\,\bullet \atop{\displaystyle \displaystyle\Downarrow \atop\displaystyle\bullet}}}

\def\mathcalB{  \star \lra \star\leftrightarrow\star  }
\def\mathcalB{\left(\displaystyle{\overset{\star}{\displaystyle\Downarrow}\,\,\,\,\,\,\,\,}\atop{ \star\leftrightarrow\star }\right)}
\def\mathcalC{\left({\emptyset\atop{\displaystyle\Downarrow\atop\displaystyle  \bullet }}\right)}
\def\mathcalCnobrackets{{\emptyset\atop{\displaystyle\Downarrow\atop \displaystyle \bullet }}}

\def\mathcalE{( \bullet\leftrightarrow\bullet \lra  \bullet )}
\def\mathcalE{\left( {\star\leftrightarrow\star}\atop{\displaystyle\displaystyle\Downarrow \atop{ \star }}\right)}
\def\mathcalEnobrackets{{ \star\leftrightarrow\star \atop{\displaystyle\displaystyle\Downarrow\atop \star }}}

\def\mathcalFbullet{\filledstar}\def\mathcalFsmallblackbox{\filledstar}
\def\mathcalF{{\,\,{\displaystyle\mathcalFbullet\leftrightarrow\mathcalFbullet_\downarrow\phantom{\mathcalFsmallblackbox\leftrightarrow\mathcalFsmallblackbox}}\atop{\displaystyle\phantom{\star\star}\displaystyle\mathcalFsmallblackbox\leftrightarrow\mathcalFsmallblackbox}}
	     \atop{\displaystyle\Downarrow\atop{\displaystyle\star\leftrightarrow\star\leftrightarrow\star}}}

\def\mathcalFwithIndices{{\!\!{\displaystyle\mathcalFbullet_u\leftrightarrow\mathcalFbullet_{\!\!v\atop\downarrow}\phantom{\mathcalFsmallblackbox\leftrightarrow\mathcalFsmallblackbox}}\atop{\displaystyle\phantom{\star_u\leftrightarrow\star_v}\displaystyle\mathcalFsmallblackbox_c\leftrightarrow\mathcalFsmallblackbox_d}}
	     \atop{\displaystyle\Downarrow\atop{\displaystyle\star_u\leftrightarrow\star_{v=c}\leftrightarrow\star_d}}}

\def\mathcalFbrackets{\left({\,\,{\displaystyle\mathcalFbullet\leftrightarrow\mathcalFbullet_\downarrow\phantom{\mathcalFsmallblackbox\leftrightarrow\mathcalFsmallblackbox}}\atop{\displaystyle\phantom{\star\star}\displaystyle\mathcalFsmallblackbox\leftrightarrow\mathcalFsmallblackbox}}
	     \atop{\displaystyle\Downarrow\atop{\displaystyle\star\leftrightarrow\star\leftrightarrow\star}}\right)}

\def\mathcalFbracketsSection{\left({\!\!{\displaystyle\mathcalFbullet\leftrightarrow\mathcalFbullet_\downarrow\,\,\,\,\,\,\,}\atop{\displaystyle\,\,\,\,\,\,\,\,\,\,\,\,\,\,\,\,\,\,\mathcalFsmallblackbox\leftrightarrow\mathcalFsmallblackbox}}
	     \atop{\displaystyle\Downarrow\atop{\displaystyle\star\leftrightarrow\star\leftrightarrow\star}}\right)}

\def\mathcalG{\left({\displaystyle\bullet\phantom{,\bullet}}\atop{{\!\!\!\displaystyle\Downarrow}\atop{\displaystyle \bullet\,\,\bullet}}\right)}

\def\mathcalGnobrackets{{\displaystyle \bullet\phantom{,\bullet}}\atop{{\!\!\!\displaystyle\Downarrow}\atop{\displaystyle \bullet\,\,\bullet}}}

\def\mathcalGnobracketsSection{{\displaystyle \bullet\,\,\,\,\,\,\,\,}\atop{{\!\!\!\displaystyle\Downarrow}\atop{\displaystyle \bullet\,\,\bullet}}}

 \def\closedsubspace{{{ \raisebox{2pt}{$\star_x\leftrightarrow\star_y\leftrightarrow\star_z$} \searrow \raisebox{-2pt}{$\smallblackbox_c$} }\atop{\displaystyle\Downarrow}}
 \atop{  \star_x=\star_y\leftrightarrow \star_z=\bblackbox_c }}

 \def\closedsubspace{{{\raisebox{2pt}{$\star_x\leftrightarrow\star_y\leftrightarrow\star_z$} \searrow \raisebox{-2pt}{$\smallblackbox_c$} }\atop{\displaystyle\Downarrow}}
 \atop{\displaystyle  \star_{x=y}\leftrightarrow \star_{z=c} }}

\def\mathcalIso{{\emptyset\atop{\displaystyle\Downarrow\atop \displaystyle \emptyset }}}

\def\QN{\mathcal{QN}}

\begin{document}

\begin{abstract}
%
The free monoid with two generators acts on classes (=properties) of morphisms of a category by 
taking the left or right orthogonal complement with respect to the lifting property, 
and we define the {\em Quillen negation monoid} of the category to be 
its largest quotient which acts 
faithfully. 
We consider the category of topological spaces and show
that a number of 
natural properties of continuous maps are obtained by 
applying this action to a single example.

Namely, for the category of topological spaces 
we show finiteness of the orbit of the simplest class of morphisms 
$\{\emptyset\lra\{\star\}\}$, 
and we calculate its Schreier graph.

 The orbit consists of 21 classes of morphisms, 
 and most of these classes are explicitly defined by standard terminology
 from a typical first year course of topology: 
 a map having {a section} or {dense image}; 
 quotient and induced topology; 
 surjective, injective; 
 (maps representing) 
 subsets, closed subsets; 
 disjoint union, disjoint union with a discrete space; 
 each fibre satisfying separation axiom $T_0$ or $T_1$. 
 Also, the notions of being connected, having a generic point, and 
 being a complete lattice, can be defined in terms of the classes 
 in the orbit.

In particular, calculating parts of this orbit can be used in an introductory course 
as  exercises connecting basic definitions in topology and category theory.

\end{abstract}
\maketitle

\begin{flushright}\tiny
\begin{minipage}[t]{0.73\textwidth}
\begin{minipage}[t]{0.958\textwidth}
{\tiny	Die Mathematiker sind eine Art Franzosen: Redet man zu ihnen, so
\"ubersetzen sie es in ihre Sprache, und dann ist es alsobald ganz etwas
	anderes.\ } 
\\	{\tiny --- Johann Wolfgang von Goethe. Aphorismen und Aufzeichnungen. Nach den Handschriften des Goethe- und Schiller-Archivs hg.~von Max Hecker, Verlag der Goethe-Gesellschaft,
  Weimar 1907.
  Aus dem Nachlass, Nr.~1005, Uber Natur und Naturwissenschaft. Maximen und Reflexionen.
}%
\end{minipage}\end{minipage}\end{flushright}

\section{Introduction}
We say that {a property (=class) $P$} of morphisms in a category $\mathcal C$ is {\em Quillen definable}
in terms of property $Q$ iff the class $P$ can be obtained from the class $Q$ by repeatedly
taking the left and right orthogonal completement with respect to lifting property, 
a binary relation on morphisms of a category used in a prominent way by Quillen in 
an axiomatic approach to homotopy theory  
\cite[1\S1, Def.1; 2\S2($\text{sSets}$), Def.1, p.2.2; 3\S3($\Topp$), Lemma 1,2, p.3.2]{Qui}.


A number of standard basic properties of morphisms 
are Quillen-definable 
in terms of (the class consisting of) a single and simple example of a morphism
violating or satisfying the property. And so are standard basic properties of objects,
if, somewhat informally, we say that a property $P'$ of objects is Quillen definable in terms of $Q$
iff it is of form $f_X \in P$ where $P$ is Quillen definable in terms of $Q$ and 
$f_X$ denotes a morphism associated with an object $X$ in some way, e.g. 
the initial or terminal morphism with (co)domain $X$. 

In this way one can define the properties of 
a (finite) group being nilpotent, soluble, $p$-group \cite[Corollary 3.2]{CR}; 
a module being injective or projective; 
a metric space being complete;
a topological space being compact, contractible, connected, having a generic point, 
totally or extremally disconnected, (hereditary) normal, (in)discrete \cite{vita,YetAnotherNotAgain}.
Moreover, for the topological properties listed above 
the simple examples can 
very simple indeed: 
maps of finite topological spaces of size $\leqslant 5$. 
This leads to a concise combinatorial notation for properties of continuous maps 
allowing to encode 
quite a few of standard topological definitions
in 2 or 4 bytes.  

In this paper we consider the category of topological spaces and 
classify the properties Quillen-definable 
in terms of the (class consisting of) the single morphism 
$ \emptyset \lra \{ \star \}  $;  there are precisely 21 of them 
and most are explicitly introduced in a first year course of topology, 
such as a map having a dense image, having a section, being a (closed) immersion.
More explicitly, 
we calculate the orbit of  $\{ \emptyset \lra \{ \star \} \} $ 
and its Schreier graph (Fig.~\ref{Fig:OrbitWithLabels}) 
under the action of the free monoid with two generators $l$ and $r$
taking a class of morphisms into its left or right orthogonal complement
$$P \longmapsto P^l \ \ \ \ \ P \longmapsto P^r$$  
We define the {\em Quillen negation monoid} of a category to 
be the largest quotient such that this action is faithful;
a justification for this definition there are interesting examples
of orbits of its action.
We formulate a couple of questions, e.g. 
whether the Quillen negation monoid is finite for 
the category of topological spaces or for the category
of finite groups of a fixed period. 

Thus one may ``stumble upon'' or ``generate'' these notions with a non-negligible probability
simply by picking an interesting example of a map or a pariculalry simple class of maps, 
and applying the trick a few times. This suggests that Quillen negation can be viewed
as a rule of ergologic of Gromov \cite{Ergobrain}, see \cite{DMG,YetAnotherNotAgain} for detailed speculations.

Let us now repeat the above in more detail.

\subsubsection*{A sketch of our results} In this note 
we look at the category of topological spaces
 and the notions defined using the lifting property in terms of the simplest morphism
 $\emptyset\to \{\bullet\}$, 
 namely the map from the empty set to the singleton. We find that most of them  are quite natural 
 and there are only finitely many of them, in a precise sense we describe now.

Following \cite[1\S5,Def.1(M6)]{Qui}, see Def.~\ref{def:rtt} below, for a 
 property (=class) $P$ of morphisms in a category $\mathcal C$
 we define the properties (=classes) $P^l$ and $P^r$ 
 of having the left, resp.~right, 
 lifting property with respect to each morphism in $P$. 
 This defines an action of the free monoid with two generators $l$ and $r$ 
 on properties (=classes) of morphisms in the category. 
We say that {\em a property $P$ Quillen defines a property $Q$} iff $Q$ lies in the orbit of $P$, i.e.
$Q=P^w$ for some word $w\in \{l,r\}^{{\aleph_0}}$. 
We also say that {\em property $P$ $w$-defines property $Q$} when $Q=P^w$. 
Often we consider $P=\{f\}$ consisting of a single map of finite spaces, and then 
we say that {\em a map $f$ $w$-defines a property (=class) $Q$ of morphisms} 
iff $Q=\{f\}^{w} $. 
We also say that {\em a map $f$ $w$-defines a property  $Q$ of spaces for spaces with property $R$}
iff for each space $X$ with property $R$ it holds ${ X\atop{\downarrow\atop\ptt} }\in \{f\}^{w} $ iff $X$ has property $Q$.

In this paper we classify the properties defined by the map $\emptyset\lra\ptt$, or, equivalently,
calculate the orbit of $P=\{\emptyset\lra\ptt\}$.
 Our main result is Theorem~\ref{thm:FinOrb}
 saying that  its orbit is finite, and consists of 21 properties. Fig.~\ref{Fig:Orbit} and Fig.~\ref{Fig:OrbitWithLabels} show
 the Schreier graph of the orbit. 
Moreover, 7 of these properties are explicitly defined by standard terminology introduced 
in a typical first year 
course of topology 
(namely,  quotient, induced topology, subspace, closed subspace, 
having dense image, having a section, surjective, injective, disjoint union),
and 5 more almost so (disjoint union with a discrete space, having a section picking a generic point in each fibre, 
each fibre satisfies Separation Axiom $T_0$ or $T_1$, the domain is (non)-empty), and 
the map $\pi_0(f): \pi_0(X)\lra \pi_0(Y)$ is surjective (only for those $Y$ where $\pi_0(Y)$ is finite)). 
As each element of the orbit can be coded by a $l,r$-word
of length at most $7$, we see that this action leads to a notation so concise that 
each of these notions fits into a single byte. 

The main tool is the proof is a concise combinatorial notation for maps of finite topological spaces
using the fact that a finite topological space is the same as a finite preorder.
In particular, $8$ of the properties in the orbit are of form $\{f\}^l$, and $4$ of form $\{f\}^r$,
where $f$ is a map of finite spaces of size $\leqslant  4$,
and in our calculation 
we specify these maps explicitly  in \S\ref{sec:negations_of_maps}. 
In fact, we believe  it should not be hard to extend our result and and calculate explicitly 
the likely finite Schreier graph 
generated by the classes of maps of spaces with $\leq 2$ elements. 
Extending this to maps of spaces of size $\geq 3$ might be more difficult, as the following examples show.
For Hausdorff spaces compactness is $lr$-defined by a map from $3$ point space to $2$ points, 
and contractibility for finite CW complexes is $lr$-defined by a class consisting of two maps from $\leq 5$ points to $\leq 3$ points. 
(\cite[Lemma 2.3.1]{vita}, \cite[Theorem~\ref{coro:contrCW}, Corollary~\ref{coro:KtoOlrSimplifiedFurther}]{YetAnotherNotAgain}).

\subsubsection*{Conclusion.} Together with the action described above,
our notation for monotone maps of finite preorders
leads to a concise combinatorial notation for the basic topological notions mentioned above, but also
for topological properties such as compact, contractible, connected, $\pi_0(f)$ being surjective or injective, 
and, conjecturally, 
a proper map, a trivial fibration 
(among maps of ``nice'' spaces).

This shows that there is combinatorics of finite preorders
implicit in  basic definitions of topology. We explore this further in \cite{vita,YetAnotherNotAgain}.

\subsubsection*{Structure of the paper} In \S\ref{intro:QM} we define the action of a monoid on a category,
and in \S\ref{intro:Figs} describe the orbit of $\{\emptyset\lra\ptt\}$ in $\Topp$. 
In \S\ref{intro:OpenProblems} we formulate several open problems.
In \S\ref{prelim:rtt} we define the lifting property, and in \S\ref{prelim:FinTop} introduce the notation for maps of finite topological spaces.
In \S\ref{prelim:qnSets} we demonstrate our ideas and our notation describing  the ``baby'' example of the Quillen negation monoid of the category of Sets.
Then in \S\ref{sec:negations_of_maps} we rewrite several properties 
in terms of maps of finite topological spaces (=preorders) of size $\leqslant 4$,
and use them in to calculate the iterated Quillen negations/orthogonals of $\{\emptyset\lra\{\bullet\}\}$ in \S\ref{sec:lrls}. 

\subsubsection*{Acknowledgements} We thank Anna Erscher for discussions and remarks 
on a preliminary version of the paper, in particular for pointing out 
that our results are best formulated in terms of the Quillen negation monoid and its action. 
We are grateful to Kobi Peterzil for support and hospitality.


\subsection{The Quillen negation monoid of a category\label{intro:QM}}
We introduce a monoid associated with a category, and 
formulate our results in terms of its action on classes of morphisms in the category.
As we explain below, this action is generated by two transformations
$P\longmapsto P^l$ and $P\longmapsto P^r$ sending a class of morphisms
into its left or right orthogonal. 

\subsubsection{The definition of a Quillen negation monoid of a category}
Fix a category $\mathcal C$. 
 A simple way to define a class of morphisms without a given property $P$
 is by taking its {\em left} or {\em right orthogonal} or $P^l$ or $P^r$ with respect 
 to the lifting property (orthogonality of morphisms, see Def.~\ref{def:rtt}), defined 
 as the class of all morphisms which have the left, respectively right, lifting property with respect to each morphism with property $P$.
It is convenient to refer to (the property of being in the class) $P^l$ and $P^r$ as
the {\em left} or {\em right} {\em Quillen negation} of (the property of being in the class) $P$.
Taking the orthogonal/Quillen negation $P\mapsto P^l$ and $P\mapsto P^r$ defines
 the action of the free monoid with 2 generators $l$ and $r$ on 
 classes (=properties) of morphisms in a category: a word $w$, e.g. $w=lrrr$,
 sends a property $P$ into $P^w$, e.g.~$P^{lrrr}:=(((P^l)^r)^r)^r$;
 in this situation it is convenient to say that {\em property $P$ $w$-defines property $Q=P^w$}.

The {\em Quillen negation monoid} $\QN(\mathcal C)$ is the quotient monoid making this action faithful:
\begin{defi}
	The  {\em Quillen negation monoid} $\QN(\mathcal C)$ of a category $\mathcal C$ consists of words in letters $l$ and $r$.
Two words are considered equal iff they act the same on each property, i.e. 
$$
v = _{\QN(\mathcal C)} w \iff \text{ for each property }P\text{ } P^{v}=P^{w}$$

Fix a property $P_0$ of morphisms in the category $\mathcal C$. 
Two words are equal in the  {\em Quillen negation monoid
	$\QN(\mathcal C,P_0)$ of a property $P_0$} 
iff they act the same on each property in the orbit of $P_0$, i.e.	$$
v = _{\QN(\mathcal C,P)} w \iff \text{ for each word }s\text{ } P_0^{sv}=P_0^{sw}$$
\end{defi}
An easy calculation (see \S\ref{prelim:qnSets} for a picture of the Schreier graph) shows that the Quillen negation monoid $\QN(\Sets)$ of the category of sets 
is finite, and that there are precisely 8 different properties of maps in $\Sets$ defined by Quillen negation (i.e.~of form $P^l$ or $P^r$
for some class $P$); 
essentially it is the same calculation which shows that there are precisely 9 model structures on the category of $\Sets$ \cite{GW}.

\subsubsection{Definitions in terms of Quillen negation} 
A number of standard basic definitions 
can be concisely and uniformly expressed 
as $ P^w$ or $f_G\in P^w$ 
where $w\in \QN(C)$, 
the morphism $f_G$ is constructed explicitly from an object $G$, 
and $P$ is a class of simple examples
of morphisms satisfying or violating a property related to the definition.
This is how Quillen defined (co)fibrations and acyclic (co)fibrations (for $w=l,lr$) in several model categories
in his axiomatic approach to homotopy theory  \cite[2\S2($\text{sSets}$), Def.1, p.2.2; 3\S3($\Topp$), Lemma 1,2, p.3.2]{Qui}.
In the category of (finite) groups, 
examples include   the definitions of 
   a finite group being nilpotent, solvable,
      perfect, torsion-free; $p$-groups, and prime-to-$p$ groups, perfect core, Fitting subgroup, and $p$-core \cite{CR}.
In model theory, 	    
Shelah's characterisations of stability, NIP, NOP, and non-dividing are of this form \cite{Z1,Z2,S}.
%
\subsubsection{Topological properties defined by iterated Quillen negation} 
In the category of topological spaces it is often enough to take $P$ to be a finite class of maps of finite topological spaces, say of size at most 6.
This leads to a concise combinatorial notation for properties such as 
a topological space being compact, contractible, connected, zero-dimensional, 
separation Axioms $T_0$, $T_1$, $T_4$ (normal), and $T_6$ (hereditary normal); a map having dense image, being a closed inclusion, and 
having connected fibres, being proper (for maps of nice enough spaces, e.g. metrisable) \cite[\S1.1]{YetAnotherNotAgain},\cite{LP2}.

\subsection{The finite orbit of $\{\emptyset\lra \ptt\}$ in the category of topological spaces\label{intro:Figs}}
In this paper we calculate the orbit of the action of the Quillen negation monoid
on the simplest class, namely the class $P:=\{\emptyset\to \{\bullet\}\}$ consisting 
 of a single morphism $\emptyset\to \{\bullet\}$ from the empty set to the singleton.
 It turns out this orbit is finite and consists of 21 rather natural properties of maps,
8 of which are explicitly defined in a typical introductory course of topology, and about 6 more almost so. 
For example, the word $rllrrl$ sends  $\{\emptyset\to\{\bullet\}\}$ into the class of closed inclusions
representing the notion of a ``closed subset'', and the word $rllrrll$ sends it 
into the class of maps having the dense image; see \S\ref{orb:meanings} for a list. 

\begin{theo}\label{thm:FinOrb} The orbit of $\{\emptyset\lra\ptt\}$ under the action of the Quillen negation monoid $\QN(\Topp)$ of the category of topological spaces,
	 is finite and has 21 element. 
\end{theo}
\begin{proof} The orbit is calculated case by case in \S\ref{sec:lrls}. \end{proof}

\subsubsection{The Schreier graph of the orbit} 
\begin{figure}
\vskip-2pt\noindent\includegraphicss[1.04949]{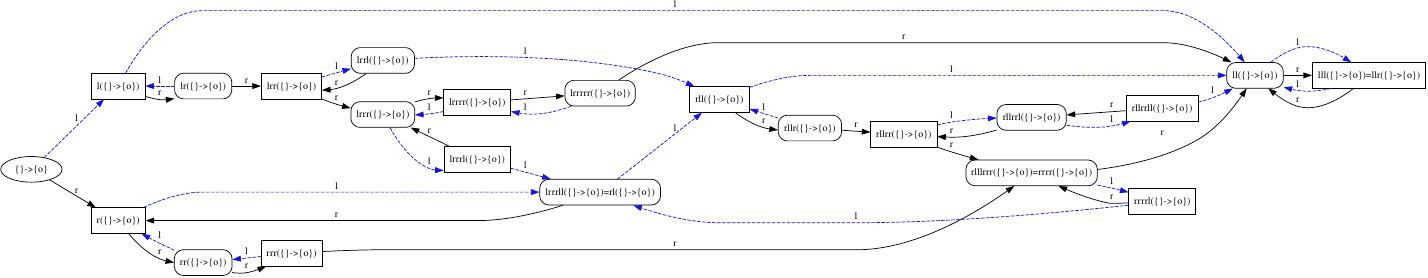}
	\caption{The orbit of $\{\emptyset\lra\ptt\}$ in $\Topp$.\label{Fig:Orbit}} 
\end{figure}

\begin{figure}
	\vskip-2pt\noindent\includegraphicss[1.0010949]{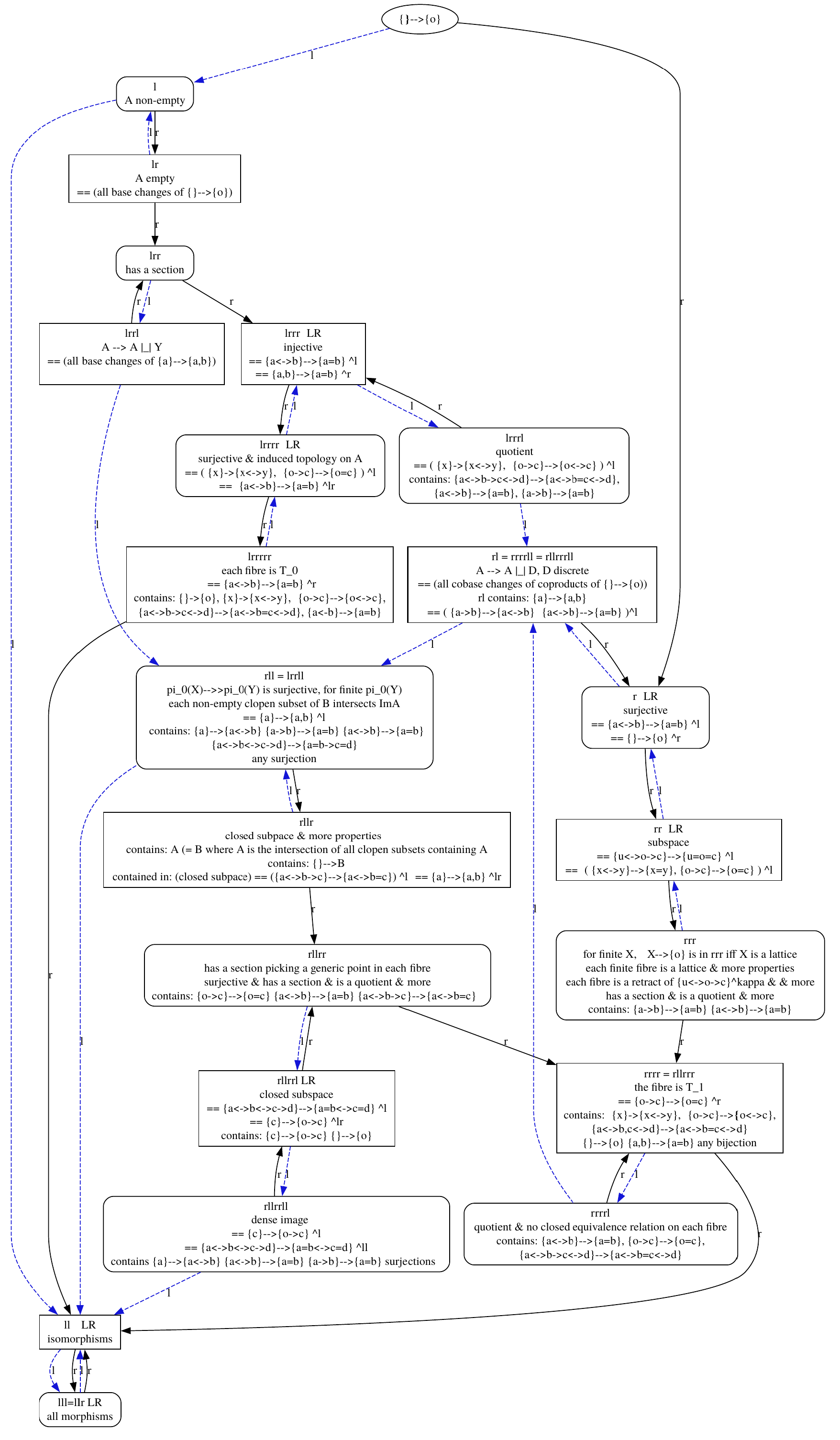} 
	\caption{The orbit of $\{\emptyset\lra\ptt\}$ in $\Topp$.\label{Fig:OrbitWithLabels}} 
\end{figure}

\begin{figure}
\tiny.\hskip-22pt\vskip-32pt\noindent\includegraphicss[1.2949]{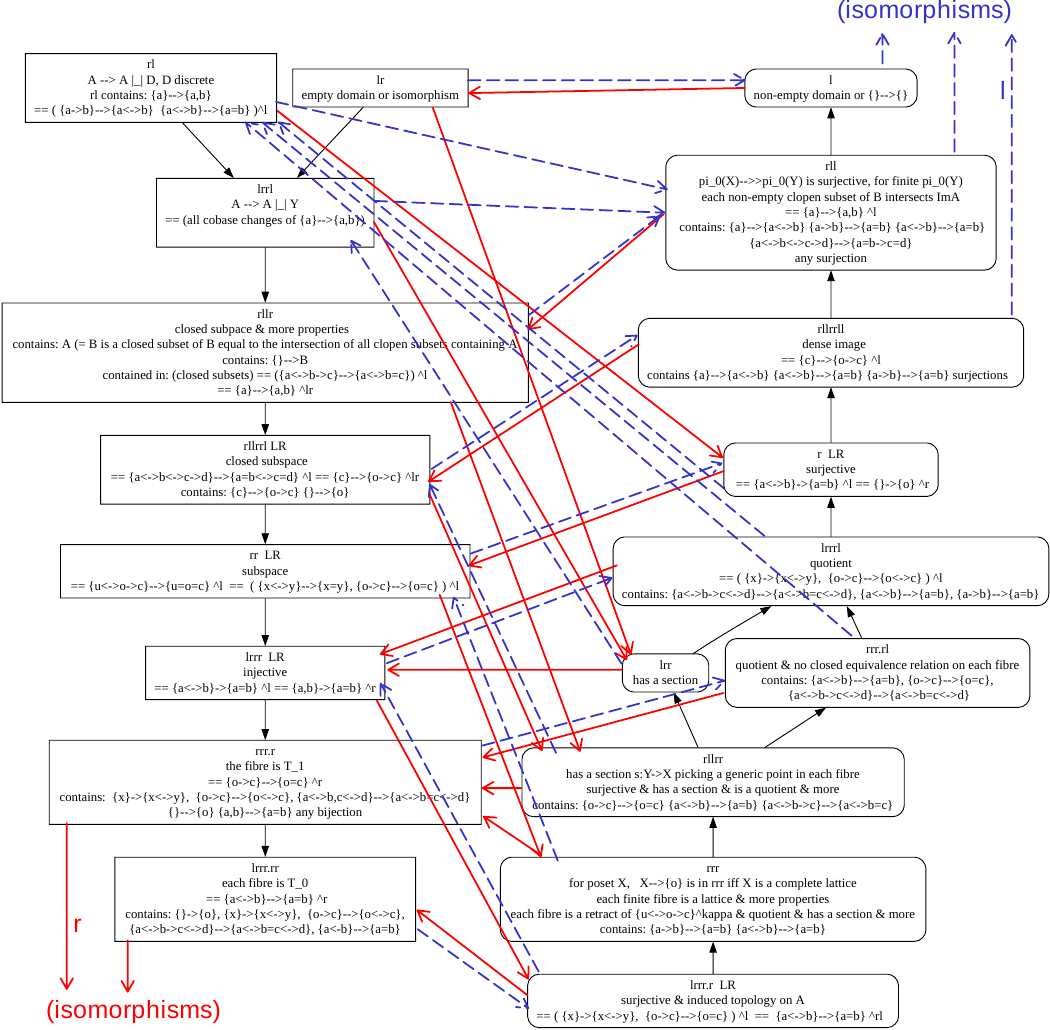}%
	\caption{The classes in the orbit of $\{\emptyset\lra\ptt\}$ in $\Topp$
	as a bipartie graph with sides ordered by inclusion. Black arrows represent inclusion, $r$-arrows are red, 
	and $l$-arrows are blue dashed.
	 \label{Fig:Orbit-Inclusions}} 
\end{figure}

\begin{figure}
\tiny.\hskip-22pt\vskip-32pt\noindent\includegraphicss[1.2949]{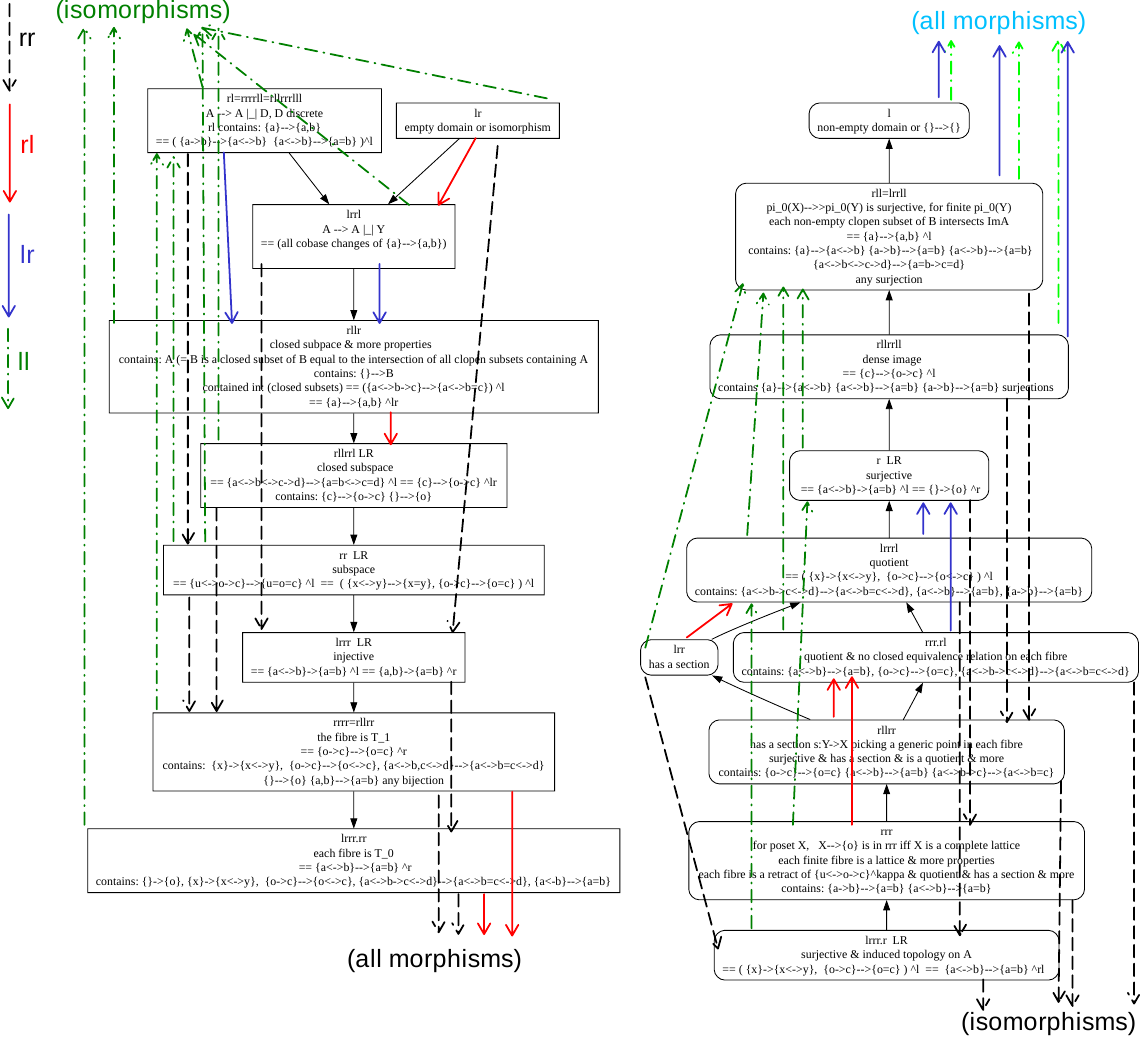}%
	\caption{The classes in the orbit of $\{\emptyset\lra\ptt\}$ with generators $ll,lr,rl,rr$ in $\Topp$
	as a bipartie graph with sides ordered by inclusion. Black arrows represent inclusion;
	loops are not depicted.
	 \label{Fig:OrbitLLRR-Inclusions}} 
\end{figure}

Fig.~\ref{Fig:Orbit}  represents the Schreier graph  of the orbit of $\{\emptyset\lra\ptt\}$ 
under $\QN(\Topp,\{\emptyset\lra\ptt\})$ action. The label on a vertex represents (usually the shortest) path from the ``root'' $\{\emptyset\lra\ptt\}$
to the vertex. Fig.~\ref{Fig:OrbitWithLabels} shows the same Schreier graph where we add to each vertex an informal explanation of the class it represents.
Fig.~\ref{Fig:Orbit-Inclusions} and Fig.~\ref{Fig:OrbitLLRR-Inclusions} 
shows that the Schreier graph is a bipartie graph with sides represented 
by words of even and odd length (if we ignore two vertices corresponding to 
the classes of isomorphisms and of all morphisms), 
and that each side is ordered by inclusion almost linearly. 
The graph has $12$ cycles of length $2$, and 
two cycles of length $6$ sharing the path $rrrr=rllrrr \to rrrrl\to rl$:
$$ll rrrr\,: \,\,\, rrrr \xrightarrow l rrrrl \xrightarrow l rl \xrightarrow r r \xrightarrow r rr   \xrightarrow r rrr   \xrightarrow r rrrr \,\,\,\,\,\,\,$$ 
$$lll rrr\,:\,\,\, rrrr \xrightarrow l rrrrl \xrightarrow l  rl \xrightarrow l rll \xrightarrow r  rllr \xrightarrow r rllrr \xrightarrow r  rrrr $$ 
Each 2-cycle is formed by the $l$- and $r$-arrows going in opposite directions. 

%
%

The pair $ll,  lll=llr$ forms a ``sink'': if a path enters either of the two vertices, 
it starts forever going back and forth $ll\leftrightarrow lll=llr$.
Each word containing $l^6$ or $r^6$ enters the sink $ll\leftrightarrow lll$, 
hence for each $s$ there are exactly two distinct words of form $sl^6...$ and of form $sr^6...$ (namely, of odd length and of even length).

The 
graph becomes a tree if we leave only $l$-edges or only $r$-edges and remove the vertex $lll$.


It is easy to describe the behaviour of a path though this graph. First, we may remove the  $lr$ 2-cycles. 
A path without $lr$ 2-cycles either ends up in $ll(isomorphisms)$ in $<7$ steps, or 
enters either of the two $6$-cycles via $rl$, $rll$, or $r$, and cycles there.
It may leave only though three arrows, two of which lead to $ll$(isomorphisms) ($rll\xrightarrow l ll$, $rrrr \xrightarrow r  ll$).
If it leaves through the remaining arrow, the path ends up in $ll$ in 3 steps ($rllrr\xrightarrow l rllrrl \xrightarrow l rllrrll \xrightarrow l ll$).

The longest cycle-free path 
has 11 vertices. \comment{A computer experiment suggests that in $\QN(\Topp,\{\emptyset\lra\{\bullet\})$ 
there is a unique word of length 29  not equal to a shorter word; and in a slighly larger monoid $\QN\left(\Topp,\left\{\, \mathcalCnobrackets, 
{\mathcalAnobrackets}, {\mathcalGnobrackets}, {\OtoOc}, {\nondenseimage} \right\}  \right)$
there is a unique word of length 30 not equal to a shorter word. 
We make a small table below to make it easier to visualise:
}\comment{
rllrrrrlllrrrllrrrllllrrrllrr      length 29 
.||....|||...||...||||...||..      QN( Top, {}-->{o} )
rllrrrrlllrrrllrrrllllrrrllrr    3 words of length 29 which cannot be shortened
lrrrrllllrrrllrrrllllrrrllrrr    QN(Top, { {}-->{o}, {a,b}-->{a=b},{a}-->{a,b}, 
llrrrrlllrrrllrrrllllrrrllrrr             {o}-->{o->c},{c}-->{o->c}} )
rllrrrrlllrrrllrrrllllrrrllrrr
.||....|||...||...||||...||.. 
 |....||||...||...||||...||... 
 ||....|||...||...||||...||... 
.||....|||...||...||||...||...        length 30
\end{verbatim}
}

%

\comment{
\begin{figure}
\tiny.\vskip-32pt\noindent\includegraphicss[1.04949]{lrl-onlyL.pdf}
	\caption{The $l$-arrows in the orbit of $\{\emptyset\lra\ptt\}$ in $\Topp$. 
	Labels on the vertices give an informal description of the classes/properties.\label{Fig:Orbit-L}} 
\end{figure}
\begin{figure}
\tiny.\vskip-2pt\noindent\includegraphicss[1.04951]{lrl-onlyR.pdf}
	\caption{The $r$-arrows in orbit of $\{\emptyset\lra\ptt\}$ in $\Topp$.  
	Labels on the vertices give an informal description of the classes/properties.\label{Fig:Orbit-R}} 
\end{figure}
}

\subsubsection{The meaning of the classes in the orbit.\label{subsubsec:orb}\label{orb:meanings}} 
Several of the classes in the orbit are Quillen negations
of maps of finite topological spaces of size at most 6.  Fig.~\ref{Fig:Table} in \S\ref{sec:negations_of_maps} below
lists these maps and indicates the properties defined by their left and right Quillen negation.
In fact, it appears easy to extend our calculation to also compute the orbits of a few other morphisms,
namely (see \S\ref{subsubsec:visuals}  for a definition of the notation) 
$$\mathcalAnobrackets\,\,\,\,\,\,{\mathcalGnobrackets}\,\,\,\,\,\,
{\closedsubspace}\,\,\,\,\,
\mathcalEnobrackets\,\,\,\,\,\, \mathcalD \,\,\,\,\,\,\mathcalDEmorphism$$
for the last three 
both their left and right Quillen negation belong to the orbit. 
\comment{We give a partial picture in Fig.~\ref{Fig:qnTopExtended}. } 
However, for reasons of space this is not done in this paper. 

Let us now give a list of words defining basic properties of maps or spaces. Complete statements can be found in \S\ref{sec:lrls}. 
\begin{itemize}\item
    $rr$ --- subspace; $rllrrl$ --- closed subspace;
    $rllrrll$ --- having dense image; 
\item  $lrrl$ --- disjoint union; $rl$ --- disjoint union with a discrete space; 
\item $rll$ defines the class of maps $f:X\lra Y$ 
    such that $\pi_0(f):\pi_0(X)\lra \pi_0(Y)$ is surjective (whenever $\pi_0(Y)$ is discrete).

\item     $lrrrl$ ---  quotient (i.e. the class of maps $f:X\lra Y$ such that the topology on $Y$ is the quotient topology);
$lrrrr$ --- the map is surjective and the topology on its domain is induced from its image.
    \item 
$rrrr$ --- each fibre satisfies Separation Axiom $T_1$; and $lrrrrr$ --- each fibre satisfies Separation Axiom $T_0$;
 
 \item 
    $lr$ --- the domain is empty; $l$ --- the domain is non-empty; 
\item 
    $r$ --- surjective; $lrrr$ --- injective; 

\item $lrr$ --- having a section; $rllrr$ --- having a section picking a generic point in each fibre.

\item $rrr$ --- a partial order $P$ is a complete lattice iff $P\lra\ptt$ lies in $\{\emptyset\to \ptt\}^{rrr}$.
\end{itemize} 

\subsection{Open problems\label{intro:OpenProblems}} 
We formulate a couple of obvious questions one may ask about the Quillen negation monoid.

\subsubsection{Randomly generating a definition of compactness and contractibility} Our observations lead to a notation for topological properties of maps or spaces
so concise that the definitions of compactness, contractility, and connectedness fit into two or four bytes.
This notation makes explicit finite preorders implicit in these notions. Let us explain.

Fix a random distribution on finite preorders. 
The following problem asks what is the probability that a map of finite spaces defines connectedness, compactness, and contractibility. 
\begin{enonce}{Problem} What is the probability that a map of finite spaces of size $\leqslant n$ 
	``defines'' connectedness, compactness, or contractibility, in the sense that it satisfies either of the conditions (i),(ii), and (iii), resp. 
\begin{itemize} 	\item[(i)] A space $X$ is connected iff $X\lra\ptt \,\in\, \{g\}^l$, i.e. $X\lra \ptt\,\rtt\, g$.
\item[(ii)] A Hausdorff space $X$  is compact iff $X\lra \ptt\,\in \{g\}^{lr}$. 

\item	[(iii)] A finite CW complex $X$  is contractible  iff $X\lra \ptt\,\in \{g\}^{lr}$. 
\end{itemize}
	What is the probability that $\{f\}^l$ or $\{f\}^{lr}$ lies in the orbit of $\{\emptyset\lra\ptt\}$ ? 
\end{enonce} 
\begin{proof}[Explanation]
	A verification shows that a map satisfies (i) iff each of its fibres is discrete, and 
	at least with one of fibres has at least two points.\footnote{Indeed, 
	in the lifting property $X\lra \ptt \,\rtt\, g$ we only need to consider
	        the fibres of $g$: in the commutative square 
	the space $X$ has to map to a fibre of $g$, and thus it holds iff for each fibre $F$ of $g$ it holds $X\lra Y \,\rtt\, F\lra\{\bullet\}$. 
	This fails if $X$ is a connected subset of $F$, which exists iff $F$ is not discrete.
	 \cite[I\S11.2,Proposition 5]{Bourbaki} implies for $F$ discrete this lifting property defines connectedness.} 

	By \cite[Corollary ~\ref{coro:proper}, Corollary~\ref{coro:KtoOlr}]{YetAnotherNotAgain} gives a purely combinatorial condition implying (ii). 
Namely, they say that	        compactness is defined as left-then-right Quillen negation of any closed (=proper) map of finite topological spaces complicated enough;
			``complicated enough'' here means that the map has as retracts the maps in \cite[Eq.~\ref{KtoOlrSimplified}]{YetAnotherNotAgain}. Similarly,
			        contractibility is defined by any trivial Serre fibration complicated enough in a similar sense using the maps mentioned in 
				\cite[Theorem~\ref{thm:Thm2}]{YetAnotherNotAgain}.
				Presumably the proportion of such maps is non-negligible.
\end{proof}
\subsubsection{A concise notation for topological properties} 
A concise and in some way intuitive notation for basic topological spaces is provided by iterated Quillen negations/orthogonals of maps of finite spaces:
{\em a word in two letters $l$,$r$, and a set of maps of finite topological spaces} represents a property (=class) of continuous maps,
By considering the morphism $X\lra\{\bullet\}$ this notation also defines a property of spaces.
\cite[Prop.~\ref{propConnDisconn},Thm.~\ref{coro:contrCW},Cor.~\ref{coro:proper}]{YetAnotherNotAgain} shows 
that to define connectedness, contractibility, and compactness, it is enough to consider one or two maps
of finite spaces of size $\leqslant 5$ and $\leqslant 3$.
\cite{LP1} gives a list of 20 topological properties defined in this way
using a single map of spaces with $\leqslant 4$ and $\leqslant 3$ points, and this paper lists some 10 properties defined
starting with the single map $\emptyset\lra\ptt$ using up to 7 Quillen negations.
In particular, \cite[Prop.~\ref{propConnDisconn}, Cor.~\ref{coro:proper}]{YetAnotherNotAgain} shows connectedness can be defined using a single map of spaces with two points,
and compactness using a single map of spaces with $4$ and $2$ points. 
A rough count on the number of maps of preorders
suggests that the definitions of these notions fit into two bytes, or perhaps three.\footnote
{Let us bound the number $\#Maps_{p\to q}$  of maps from a preorder with $p$ elements to a preoder with $q$ elements
in terms of the number of labelled preorders with $p$ elements.

Pick a preorder $Q$ with $l$ elements labelled by $1,..l$.
Pick a $P$ be a preorder labelled by $1,..,p$. A subset (=increasing sequence)  $1\leq i_1\leq ... \leq i_{l'}\leq p$ with $l'\leqslant l-1$
elements determines a (possibly not monotone) map $P\lra Q$ from $P$ into $Q$. 
Each map of unlabelled preorders with $p$ and $l$ elements can be constructed in this way.
Hence, the number of maps from a preorder with $p$ elements to a preorder with $l$ elements at is most
the product of the number of labelled preorders with $p$ elements, the number of partitions of $p$ into $\leqslant l$ intervals,
and the number of (in fact, unlabelled) preorders with $l$ elements.
Using the OEIS library (sequences 
\href{https://oeis.org/A001930}{A001930} and
and \href{https://oeis.org/A000798}{A000798})
for $p=4$ and $l=2$ we get $\leqslant 355*4*3=4260\leqslant 2^{13}$, and for $p=5$ and $l=3$ we get
$\leqslant 6942* 15*9=937170\approx 1000000\leqslant 2^{20}$
Thus, if we include the $lr$-suffix, the notions of connectedness and compactness fit into 2 bytes, and contractibility may fit into 3 or 4 bytes.
%
%
}
%
It is tempting to develop a computer algebra system using an extension  of our notation. The following is an example of a concrete goal. 
\begin{enonce}{Problem}\label{rl_rrrrll}  Develop a concise combinatorial notation for topological properties, and a computer algebra system, which can state and prove that 
	$$\mathcalC^{rl}=\mathcalC^{rrrrll}=\mathcalC^{rllrrrll}$$
	$$\mathcalC^{rllrrll}= \left(\closedsubspace\right)^{ll}=\nondenseimagebrackets^l$$ 
\end{enonce}
\begin{proof}[Explanation] We prove both identities in our calculation of the orbit of $\mathcalC$. 
	The first identity is evident in Fig.~\ref{Fig:Orbit} and Fig.~\ref{Fig:OrbitWithLabels}. 
 Theorem~\ref{crllrrll} says that classes mentioned in the second identity are equal to the class of maps with dense image.
\end{proof}

\subsubsection{A cognitive experiment} Arguably, there are finite preorders implicit in basic notions of topology. Are there preorders processed by our brain \cite[\S3]{DMG} ?

The following little experiment is perhaps both feasible and informative, once we are able to reformulate in purely combinatorial terms 
a topological argument complicated enough. 

\begin{enonce}{Experiment} Explain to a 2nd year student a combinatorial construction representing a well-known topological definition or argument. 
	What associations will the student make ? Importantly, the student should not be preconditioned to think about topology.
\end{enonce}

\subsubsection{Finitetess of orbits of maps of finite spaces?} 

\begin{enonce}{Problem} Is it true that in $\QN(\Topp)$ the orbit of any map of finite topological spaces is finite ? 

In other words, is it true that for each map $f:X\lra Y$ of finite topological spaces there are only finitely many different classes of form 
	$\{f\}^w$ where $w$ is a word in alphabet $\{l,r\}$ ? 
\end{enonce}
In particular, is this true 
for the trivial Serre fibration $\MtoL$\comment{ defined in \S\ref{def:MtoL}},
and the closed map defining compactness in \cite[\S1.1.2-3]{YetAnotherNotAgain}.

\subsubsection{Is the Quillen negation monoid of $\Topp$ finite ?} 
In any category and for any property $P$ it holds $P^{lrl}=P$ and $P^{rlr}=P^r$;
this is shown by the same calculation that shows that for a vector space its dual and its triple-dual coincide, i.e. 
for any vector space $V$ $V^{***}=V^*$. Does the Quillen negation monoid satisfy any other relation ? 
In fact, is it finite ? We state these questions as a problem.

\begin{enonce}{Problem}
Is the Quillen negation monoid $\QN(\Topp)$ of the category of topological spaces finite ? 

Does the Quillen negation monoid $\QN(\Topp)$  satisfy any relation in addition to $lrl=l$ and $rlr=r$ ?
\end{enonce}

\subsubsection{A model structure on $\Topp$ defined combinatorially?}

In a closed model category, 
if a map $f:X\lra Y$ is an (acyclic) fibration, then so in any map in $\{f\}^{lr}$. 
This formally follows from the fact that (acyclic) fibrations are defined by a right lifting property.
It is tempting to ask whether one can find a single map $f$ of finite topological spaces 
such that the class $\{f\}^{lr}$ is in some sense the class of {\em all} (acyclic) fibrations.
A candidate for such an example of an acyclic fibration is discussed in
\cite{vita}, see also \cite[\S3]{YetAnotherNotAgain}: it is an acyclic Serre fibration
such that for a space $X$ nice enough, e.g.~metrisable separable absolute neighbourhood retract, 
the map $X\lra\ptt$ is an acyclic fibration iff $X\lra\ptt\,\in\,\{f\}^{lr}$.

The following is a couple of precise questions. 
\begin{enonce}{Problem}
	Find two finite classes $(C)$ and $(WC)$ of maps of finite topological spaces,
	two words $v,w\in \{l,r\}^{{\aleph_0}}$, 
	and a model structure on the category of topological spaces such that 
	$(C)^v$ is its class of fibrations, and $(WC)^w$ is its class of acyclic fibrations.
\end{enonce}

For a class $P$, let $P_{{\aleph_0}}$ denote the subclass of $P$ consisting of maps between finite spaces.
Note that  \cite[Corollary~\ref{coro:proper}]{Y} defining compactness gives  an example 
of an expression of the form used in the problem below. 

\begin{enonce}{Problem}
	Find two 
	maps $f_{wc}$ and $f_c$ of finite topological spaces, and 
	four words $v_1,v_2,w_1,w_2\in \{l,r\}^{{\aleph_0}}$, such that for each map $f:X\lra Y$ of finite CW complexes $X$ and $Y$ it holds
	\begin{itemize}
		\item $f$ is an acyclic fibration iff $f\in((\{f_{wc}\}^{v_1})_{{\aleph_0}})^{v_2}$

		\item $f$ is a fibration iff $f\in((\{f_{c}\}^{w_1})_{{\aleph_0}})^{w_2}$
	\end{itemize}
\end{enonce}

\def\FinGroupsN{\text{FiniteGroups}\operatorname{mod} N} 
\subsubsection{A homotopy theory for the category of finite groups of a fixed period ?}
The category of finite groups of a given period $N$ (i.e.~finite groups $G$ such that $x^N=1$ for each $x\in G$) has all finite limits and colimits,
as implied by the positive solution of the Burnside problem \cite[Corollary 3.2]{NY}.\footnote{Indeed, 
the universality condition on the restricted Burnside group says precisely that it is the coproduct of cyclic groups with  $N$ elements in this category;
this is the statement of \cite[Corollary 3.2]{NY}. The rest of the claim follows by standard arguments, as follows.
An arbitrary colimit is a quotient of the coproduct. In the category of groups limits of diagrams of finite groups are also finite.}
 Let us denote this category by $\FinGroupsN$.

\cite{CR} gives examples of Quillen negations in the category of finite groups related to the notions of
being nilpotent,
solvable, perfect, torsion-free; p-groups and prime-to-p-groups; Fitting
subgroup, perfect core, p-core, and prime-to-p core. Arguably, calculations there suggest it may be possible to classify 
all the Quillen negations and weak factorisation systems in the category of finite groups or in the category $\FinGroupsN$ 
which has better category theoretic properties (namely, has finite limits and colimits).

\begin{enonce}{Problem} Calculate the Quillen negation monoid of $\FinGroupsN$. Is it finite ? 
\end{enonce}

For $N=2$ this is true for trivial reasons: a group of period $2$ is necessarily abelian and thus is the same as a vector space over the field with two elements. 
On the category of vector spaces (of finite dimension or not), there are only 4 Quillen negations: injections, surjections, isomorphisms, and arbitrary morphisms.

Recall that a model structure on a category consists of two weak factorisation systems with certain properties \cite[Introduction, p.0.1; I\S1,Def.1(Axiom M2),p.1.1]{Qui}.
Thus examining the classification of Quillen negations or weak factorisation systems on $\FinGroupsN$ can perhaps 
lead to a definition of a homotopy theory for finite groups of a fixed exponent. 

\begin{enonce}{Problem} Find a non-trivial model structure on  $\FinGroupsN$ such that its homotopy category says something non-trivial about finite groups. 
\end{enonce}

\section{Preliminaries}

\subsection{The lifting property\label{prelim:rtt}} We define the lifting property and briefly describe its properties and intuition. Proofs 
can be found in \cite[\S3]{H}, and more examples in \cite{LP2}.  
\vskip-6pt
\begin{minipage}[c]{0.31658980337\textwidth}
	$ \xymatrix{ A \ar[rr]^{\forall t} \ar@{->}[dd]|f && C \ar[dd]|g \\ \\ B \ar[rr]|-{\forall b} \ar@{-->}[uurr]|{{\exists d}}&& D }$
\end{minipage}
\begin{minipage}[c]{0.63842785841019663\textwidth}
	\begin{defi}\label{def:rtt} A morphism $A\xra f B$ in a category
	has {\em the left lifting property with respect to a morphism $C\xra g D$}, 
	and 	 {\em $C \xra g D$ also has the right lifting property with respect to $A\xra f B$},
		 denoted by ${\displaystyle f\rtt g }$ or ${\displaystyle A \xra f B \,\rtt\, C\xra g D }$,
	iff 
	for each map $A \xrightarrow t C$ and $B \xrightarrow b D$ such that  $f\circ b=t\circ g$,
	there exists $B \xrightarrow d C$ such that  $t=f\circ d$ and $b=d\circ g$.    \end{defi}
\end{minipage}
\vskip6pt
For a class $C$ of morphisms in a category, its {\em left orthogonal} or {\em left Quillen negation} $C^{\rtt l }$ with respect to the lifting property,
respectively its {\em right orthogonal} or {\em right Quillen negation} $C^{\rtt r}$, is the class of all morphisms which have the left, respectively right, lifting
property with respect to each morphism in the class $C$:
$$ C^{\rtt l}:= \{ f \mid \forall g\in C\, f\rtt g\}$$ 
$$ C^{\rtt l}:= \{ g \mid \forall f\in C\, f\rtt g\}$$

It is clear that $ C^{\rtt l r}\supset C$, $ C^{\rtt  rl }\supset C$, $ C^{\rtt l}=C^{\rtt lrl}$, and $ C^{\rtt  r} =C^{\rtt rlr}$, and that any map in $C\cap C^{\rtt l}$ or $C \cap C^{\rtt r}$ is an isomorphism. 
The class $C^{\rtt r}$ is always closed under retracts, pullbacks, 
products (whenever they exist in the category) and composition of morphisms, and contains all isomorphisms.
Meanwhile, $ C^{\rtt  l }$ is closed under retracts, pushouts, 
coproducts and transfinite composition (filtered colimits) of morphisms (whenever they exist in the category), and also contains all isomorphisms.

\subsubsection{Intuition: lifting property as negation} 
Taking the orthogonal (Quillen negation) of a class $C$ is a simple way to define a class of morphisms excluding non-isomorphisms from $C$,
in a way which is useful in a diagram chasing computation.  A useful intuition is to think that the property of left-lifting against a class $C$ is a kind of negation of the property of being in $C$, and that right-lifting is also a kind of negation, and for this reason it is convenient to 
refer to	{\em property} $C^{\rtt l}$, resp.~$C^{\rtt r}$, as  {\em the left}, resp.~{\em right,} {\em Quillen negation of property $C$.}

\subsubsection{Weak factorisation systems} 
Quillen negations are closely related to {\em weak factorisation systems}, 
defined as a pair $\mathcal L \,\rtt\, \mathcal R$ of classes of morphisms 
(i.e. for each $f\in \mathcal L$ and each $g\in \mathcal R$ it holds $f\rtt g$) such that
each morphism $f$ decomposes as $f=f_r\circ f_l$ where $f_l\in \mathcal L$ and $f_r\in\mathcal R$.
Given a weak factorisation system  $\mathcal L \,\rtt\, \mathcal R$, so are
$\mathcal L^{rl} \,\rtt\, \mathcal L^r$ and $\mathcal R^l \,\rtt\, \mathcal R^{lr}$,
where $\mathcal L\subset \mathcal  R^l\cap \mathcal  L^{rl}$ and $\mathcal R\subset \mathcal  L^r \cap \mathcal  R^{lr}$.

A desirable property of a Quillen negation $C^l$ or $C^r$ is that
each morphism $f$ decomposes as $f=f_l\circ f_{lr}$ where $f_l\in C^l$, resp.~$f_{lr}\in C^{lr}$, 
and $f=f_{rl}\circ f_r$ where $f_{rl}\in C^{rl}$ and $f_r\in C^r$, i.e. that
{\em $(C^l, C^{lr})$, resp.~$(C^{rl},C^r)$ is a  weak factorisation system}.
The property is often proven when $C$ is a set (rather than a class) 
by a transfinite construction called the {\em Quillen small object argument}.
We also verified by hand that this property holds for most of the examples considered in this paper. 

\subsubsection{Intuition: Defining  by examples.} 
A number of basic notions from a first year course in algebra or topology 
may be expressed using the lifting property 
starting from an explicit list 
of (counter)examples, 
i.e. 
as $C^{\rtt l}$, $C^{\rtt r}$, $C^{\rtt lr}$, $C^{\rtt ll}$, ...
where the class $C$ is an explicitly given list of morphisms.
For example, in $\Sets$ for the simplest non-surjections $\emptyset\to\ptt$ and $\ptt\to\{\bullet,\bullet\}$
the Quillen negation the class $\{\emptyset\to\ptt\}^r=\{\ptt\to\{\bullet,\bullet\}\}^l$ is the class of surjections,
whereas  
for the simplest non-injection $\{\bullet,\bullet\}\to\ptt$ 
the Quillen negation the class $\{\{\bullet,\bullet\}\to\ptt\}^l=\{\{\bullet,\bullet\}\to\ptt\}^r$ is the class of injections.
In the category $R-\operatorname{Mod}$ of modules, a module $P$ is projective iff $0\to P$ is in $\{0\to R\}^{rl}$, and 
a module $I$ is injective iff $I\to 0$ is in $\{R\to 0\}^{rr}$ \cite{LP1}. 
In the category of finite groups, a finite group $H$ is nilpotent 
iff the diagonal map $H\lra H\times H$ lies in $\{0\to G: G\text{ arbitrary }\}^{lr}$, is soluble iff $H\to 0$ lies in 
$\{ A \to 0 \mid A \text{ Abelian }\}$, and of order prime to $p$/power of $p$  iff $H\to 0$ lies in
$\{ \Bbb Z/p\Bbb Z \to 0\}^r$, resp.~$\{ \Bbb Z/p\Bbb Z \to 0 \}^{rr}$ \cite[Corollary 3.2]{CR}. 
\subsubsection{Retracts of Cartesian powers as a lifting property}
nWe shall use the following simple Lemma.
\begin{lemm}\label{cartRetract}\label{XtoOretracts}  Let $\mathcal C$ be a category with a terminal object $\top$, arbitrary small products, and small $\operatorname{Hom}$-sets.
Let $A$ be an object of $\mathcal C$. Then an object $X$ is a retract of a Cartesian power of $A$ iff 
$$X \lra \top \,\in\, \{ A\lra \top \}^{lr}$$
\end{lemm}\def\id{\operatorname{id}}
\begin{proof} $\implies$: By \cite[Lemma 3.6]{H} the orthogonals are closed under retracts and by the dual to \cite[Lemma 3.4]{H}  are closed under products.
$\impliedby$: Evidently $X \lra \prod\limits_{f:X\lra A}
A \,\rtt\, A \lra \top $ (namely, for a map $f:X\lra A$ the lifting $
\prod\limits_{f:X\lra A} A \lra A$ is the projection on the $f$-th
coordinate), therefore  $X \lra \prod\limits_{f:X\lra A} A \,\rtt\, X \lra \top
$. The required retraction is the lifting  
$ \prod\limits_{f:X\lra A} A  \lra X $ for the map $\id:X\lra X$.

The following diagrams explain the proof above.
$$\begin{gathered} \xymatrix@C=2.29cm{ X \ar[d]_{(A\lra\top)^l} \ar[r]^f & A \ar[d] \\  \prod\limits_{f:X\lra A} A   \ar[r] \ar@{-->}[ru]|{\operatorname{pr}_f} & \top }
\hphantom{............} 
\xymatrix@C=2.29cm{ X \ar[r]|{\id} \ar[d]_{(A\lra\top)^l} & X \ar[d]^{(A\lra\top)^{lr}} \\  \prod\limits_{f:X\lra A} A   \ar[r] \ar@{-->}[ru]  & \top }
 \end{gathered}$$
\end{proof}

\setlength{\epigraphwidth}{0.7\textwidth} 
\epigraph{A random continuing
	 example is an understanding of finite topological spaces, an oddball topic that can
	 lend good insight to a variety of questions but that is generally not worth developing
	 in any one case because there are standard circumlocutions that avoid it.
	}{W.P.Thurston. \href{arxiv.org/pdf/9404236}{\em On Proof and Progress in Mathematics}.}
The formal syntax introduced in \S\ref{FinTop:syntax} lies at heart of this paper. 
\subsection{Notation for finite topological spaces\label{prelim:FinTop}} 
\label{subsec:notationFP}
The formal syntax introduced in \S\ref{FinTop:syntax} lies at heart of this paper. 
We follow \cite[\S5.3.1]{mintsGE}.

\subsubsection{\label{app:top-notation} Finite topological spaces as preorders and as categories} 
A {\em topological space} comes with a {\em specialisation preorder} on its points: for
points $x,y \in X$,  $x \leqslant y$ iff $y \in \operatorname{cl} x$ ($y$ is in the {\em topological closure} of $x$).
The resulting {\em preordered set} may be regarded as a {\em category} whose
{\em objects} are the points of ${X}$ and where there is a unique {\em morphism} $x{\searrow}y$ iff $y \in cl x$.

For a {\em finite topological space} $X$, the specialisation preorder or
equivalently the corresponding category uniquely determines the space: a {\em
subset} of ${X}$ is {\em closed} iff it is
{\em downward closed}, or equivalently,
is a full subcategory such that there are no morphisms going outside the subcategory.

The monotone maps (i.e. {\em functors}) are the {\em continuous maps} for this topology.

\subsubsection{A syntax to denote finite topological spaces and their maps\label{FinTop:syntax}} 

We denote a finite topological space by a list of the arrows (morphisms) in
the corresponding category; '$\leftrightarrow $' denotes an {\em isomorphism} and '$=$' denotes the {\em identity morphism}.  An arrow between two such lists
denotes a {\em continuous map} (a functor) which sends each point to the correspondingly labelled point, but possibly turning some morphisms into identity 
morphisms, thus gluing some points.

With this notation, we may display continuous functions for instance between the {\em discrete space} on two points, the {\em Sierpinski space}, the {\em antidiscrete space} and the {\em point space} as follows (where each point is understood to be mapped to the point of the same name in the next space):
$$  \begin{array}{ccccccc}
  \{a,b\} &\longrightarrow& \{a{\rightarrow}b\} &\longrightarrow& \{a\leftrightarrow b\} &\longrightarrow& \{a=b\}

  \\ 
  \verb|{a,b}| &\longrightarrow& \verb|{a->b}| &\longrightarrow& \verb|{a<->b}| &\longrightarrow& \verb|{a=b}|



  \\
  \text{(discrete space)} &\longrightarrow& \text{(Sierpinski space)} &\longrightarrow& \text{(antidiscrete space)} &\longrightarrow& \text{(single point)}
  \end{array}
  $$
The second line represents how this syntax may be typed 
in a computer algebra system. We also employ it in the pictures of Schreier graphs.  

In $A \longrightarrow  B$, each object and each morphism in $A$ necessarily appears in ${B}$ as well; sometimes we avoid listing
the same object or morphism twice. Thus
both
$$
\{a\} \longrightarrow  \{a,b\}
  \phantom{AAA} \text{ and } \phantom{AAA}
  \{a\} \longrightarrow  \{b\}
  $$
  denote the same map from a single point to the discrete space with two points.

   Each continuous map $A\lra B$ between finite spaces may be represented in this way; in the first list
   list relations between elements of $A$, and in the second list put relations between their images.
   However, note that this notation does not allow to represent {\em endomorphisms $A\lra A$}.
   We think of this limitation
   as a feature and not a bug: in a diagram chasing computation,
   endomorphisms under transitive closure lead to infinite cycles,
   and thus our notation has better chance to define a computable fragment of topology.

\subsubsection{Various conventions on naming points and depicting arrows}
While efficient, this notation is unconventional and requires some getting used to. For this reason, sometimes we employ 
more graphic notation where our notation is moved to subscripts, so to say:
points or objects are denoted by bullets $\bullet,\bblacksquare, \star,...$ with subscripts,
and the reader may think that the subscripts indicate where a point maps to or what its preimage is. 
We try to make the shape of the bullet indicate whether the point is open, closed, or neither:
$\bullet$ stands for open points (which might also be closed),
$\bblacksquare$ stands for closed points, and $\star$ stands for points which are neither open or closed.
As is usual in depicting a preorder, we also try to place $\bullet,\bblacksquare, \star,...$ so that the arrows usually go downwards.
Importantly?, this notation makes visually apparent the shape of the preorder denoted.

Thus in this graphic notation we would write
$$
  \begin{array}{ccccccc}
	    \{\bullet_a,\bullet_b\}
	         &\longrightarrow&
		    \left\{{\ensuremath{\bullet_a}}\rotatebox{-12}{\ensuremath{\to}}  \raisebox{-2pt}{\ensuremath{\bblacksquare_b}}\right\}
		        &\longrightarrow&
			  \{\star_a\leftrightarrow \star_b\}
			      &\longrightarrow&
				\{\bullet_{a=b}\}
				 \\
				    \text{(discrete space)}
				         &\longrightarrow&
					   \text{(Sierpinski space)}
					       &\longrightarrow&
					         \text{(antidiscrete space)}
						     &\longrightarrow&
						       \text{(single point)}
						         \end{array}
							 $$

\subsubsection{Visual conventions\label{subsubsec:visuals}} 
\begin{figure}
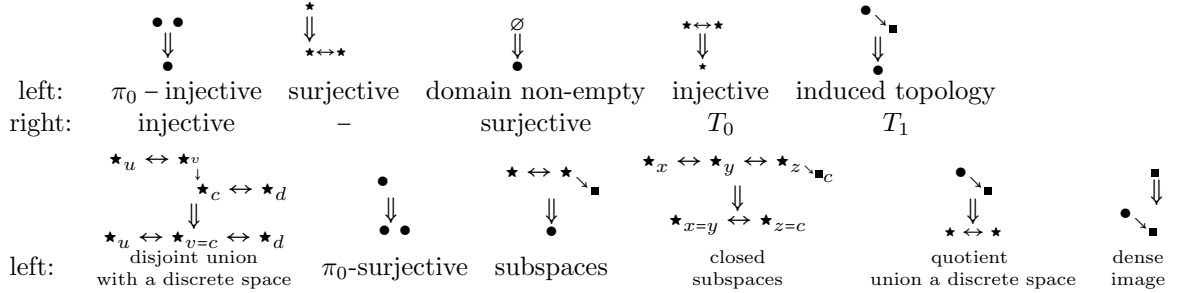

$$\begin{array}{cccccc}
&
{\mathcalAnobrackets}\hphantom{.....}&
{\mathcalBnobrackets}\hphantom{.....}&
{\mathcalCnobrackets}\hphantom{.....}&
{\mathcalEnobrackets}\hphantom{.....}&
{\mathcalD}\hphantom{.....}
\\ \text{left: } &
\pi_0-\text{injective}&
\text{surjective}&
	\text{domain non-empty}&
\text{injective}&
{\text{induced topology}}
\\ \text{right: }& 
\text{injective}&
-&
\text{surjective}&
T_0&%
T_1
\\
\end{array}
$$
$$\begin{array}{ccccccc}
&
{\mathcalFwithIndices}&
{\mathcalGnobrackets}&
{\mathcalDEmorphism}&
{\closedsubspace}&
\mathcalDD&
\nondenseimage 
\\ \text{left: } &
{\vphantom{\bar I}\text{disjoint union}\atop\text{with a discrete space}}&
\pi_0\text{-surjective}&
\text{subspaces}&%
{\text{closed}\atop{\text{subspaces}}}&
\text{quotient }\atop\text{union a discrete space}&%
{\text{dense}\atop\text{image}}
\end{array}
$$
	\caption{The Quillen negations presented as a table\label{Fig:Table}. 
	The second and third lines list the notions defined by the left, resp.~right, Quillen negation
	of these maps. See \S\ref{subsec:conn}-\ref{subsec:closed_subsets} for complete statements. }
\end{figure}

To picture a map of preorders, it also helps to place its domain above its codomain so that each point maps to a point below it.
We use $\Downarrow$ to connect the codomain and domain, as it reminds us that the preorder is a category, and a monotone map is a functor of these categories.
Whenever it does not lead to confusion, we try to skip the subscripts from this graphic notation. 
Thus, $\mathcalDD$ denotes the map from the Sierpinski space to the antidiscrete space which in the less graphic notation shall be denoted by the formula
$ \left\{{\ensuremath{\bullet_a}}\rotatebox{-12}{\ensuremath{\to}}  \raisebox{-2pt}{\ensuremath{\bblacksquare_b}}\right\} \longrightarrow 
			  \{\star_a\leftrightarrow \star_b\}$.
Fig.~\ref{Fig:Table} lists the maps we use in our calculations. The map $\mathcalAnobrackets$ glues together two points of the discrete space with two points. 
$\mathcalBnobrackets$ maps a point into the antidiscrete space with two points. 
$\mathcalD$ maps glues together the two points of the Sierpinski space. Here, in the Sierpinski space $\Otoc$
we have that $\smallblackbox\in\operatorname{cl}(\bullet)$, and  $\bullet$ is the open point, and $\smallblackbox$ is the closed point. 

In $\closedsubspace$, the two points $\filledstar_x$ and $\filledstar_y$ go to $\filledstar_{x=y}$, and 
the two points 
$\filledstar_z$ and $\smallblackbox_c$ go to $\filledstar_{z=c}$. The only non-trivial open subset of the domain 
$\raisebox{3pt}{$\filledstar_x\leftrightarrow \filledstar_y\leftrightarrow\filledstar_z$}\!\!\!\searrow \!\!\raisebox{-3pt}{$\smallblackbox_c$}$
is $\{\filledstar_x,\filledstar_y,\filledstar_z\}$. This map is neither injective nor closed. 

$\nondenseimage$ sends the singleton into the closed point of the Sierpinski space; it is an example of a map whose image is not dense, 
and also an example of a closed subset. 
\subsubsection{Conventions useful in writing up a diagram-chasing computation.}
   We
   denote points by letters $a,b,c,..,U,V,...,0,1..$ or by bullets with these subscrips 
   to make notation reflect the intended meaning,
   e.g.~an arrow $X\lra \{{{\bullet_U}\atop{}}\!\!\searrow\!\! {{}\atop{\bblacksquare_{U'}}}\!\}$ 
   reminds us that the preimage of $\bullet_U$ determines an open subset $U$ of $X$,
and   $\{\bullet_x,\bullet_y\}\lra X$ reminds us that the map determines points $x,y\in X$, and 
$\{o\searrow c\}$ reminds that $o$ is open and $c$ is closed.

\def\atob{\ensuremath{\left\{{\ensuremath{a}}\rotatebox{-12}{\ensuremath{\to}}  \raisebox{-2pt}{\ensuremath{b}}\right\}}}

\def\OtoC{\ensuremath{\left\{{\ensuremath{o}}\rotatebox{-12}{\ensuremath{\to}}  \raisebox{-2pt}{\ensuremath{c}}\right\}}}

\def\BtoSquare{\ensuremath{\left\{{\ensuremath{\bullet}}\rotatebox{-12}{\ensuremath{\to}}  \raisebox{-2pt}{\ensuremath{\bblacksquare}}\right\}}}

\def\LLambda{\ensuremath{
\left\{\underset{\bblacksquare}{}{\swarrow} \overset{\raisebox{1pt}{$ \bullet$}}{}{\searrow} 
\underset{\bblacksquare}{}\right\}}}

\def\LLambdaauxvb{\ensuremath{
\left\{\underset{\bblacksquare_a}{}{\swarrow} \overset{\raisebox{1pt}{$ \bullet_{u=x=v}$}}{}{\searrow} 
\underset{\bblacksquare_b}{}\right\}}}

\def\MMM{\ensuremath{
\left\{\underset{\bblacksquare}{}{\swarrow} \overset{\raisebox{1pt}{$ \bullet$}}{}{\searrow} 
\underset{\bblacksquare}{}
{\swarrow}\overset{\raisebox{1pt}{$\bullet$}}{}{\searrow}\underset{\bblacksquare}{}\right\}
}}

\def\MMMauxvb{\ensuremath{
\left\{\underset{\bblacksquare_a}{}{\swarrow} \overset{\raisebox{1pt}{$ \bullet_{u}$}}{}{\searrow} 
\underset{\bblacksquare_{x}}{}
{\swarrow}\overset{\raisebox{1pt}{$\bullet_{v}$}}{}{\searrow}\underset{\bblacksquare_b}{}\right\}
}}

\def\MMMABUVX{\ensuremath{
\left\{\underset{\bblacksquare_A}{}{\swarrow} \overset{\raisebox{1pt}{$ \bullet_{U}$}}{}{\searrow} 
\underset{\bblacksquare_{X}}{}
{\swarrow}\overset{\raisebox{1pt}{$\bullet_{V}$}}{}{\searrow}\underset{\bblacksquare_A}{}\right\}
}}


\subsection{The Quillen monoid of the category of Sets\label{prelim:qnSets}} 
\begin{figure}
	\hskip-12pt\noindent\includegraphicss[1.204949]{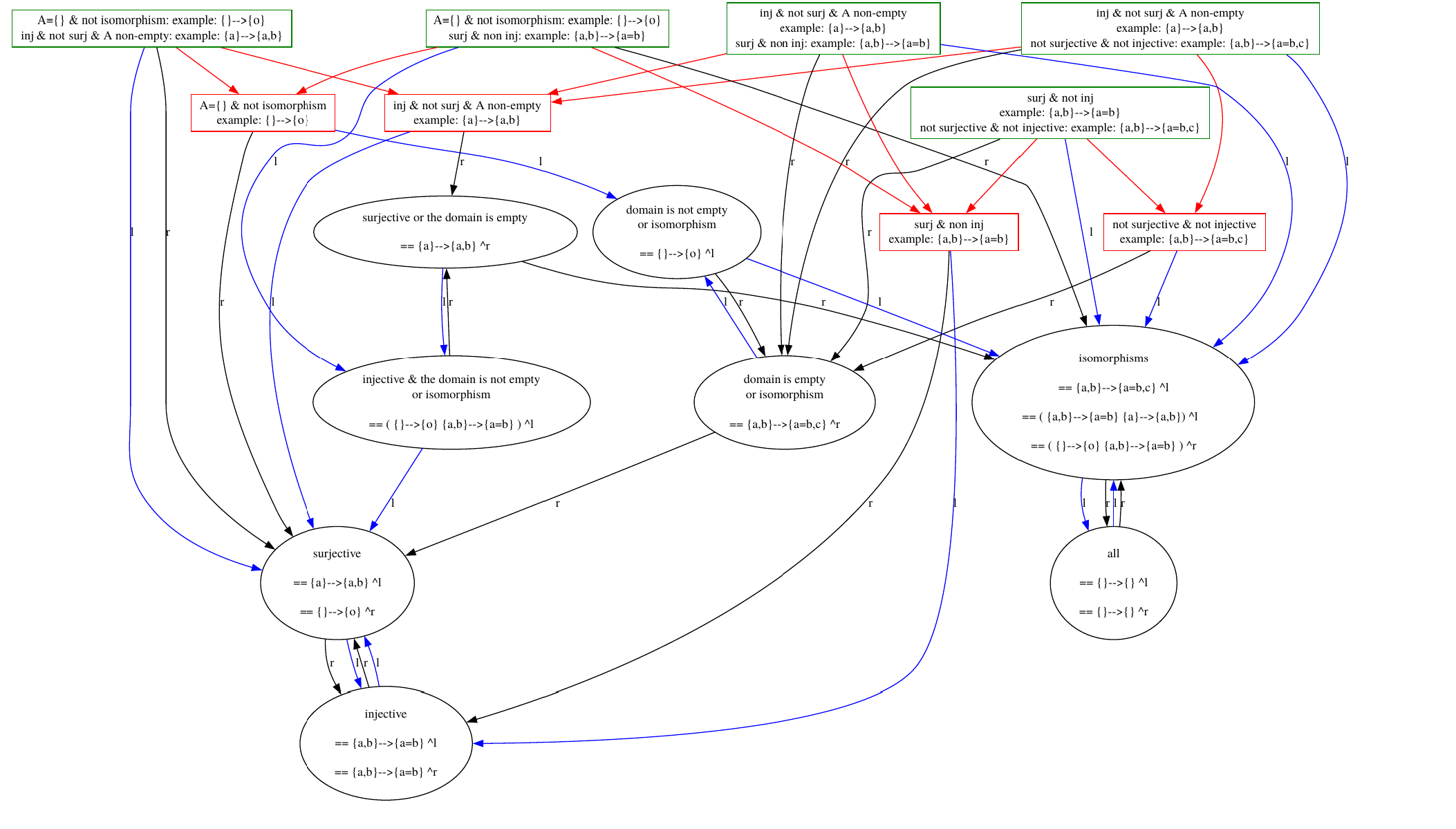}
	\caption{The Quillen negation monoid of the category of sets.\label{Fig:QnSets}} 
\end{figure}

This toy example clarifies the notion of the Quillen groupoid and, more importantly, our notation for maps of finite topological spaces.
Fig.~\ref{Fig:QnSets} represents the Quillen monoid of the category of sets. Let us now explain this picture. 

In the category of Sets there are exactly 8 classes (=properties) of maps which have form of a Quillen negation $P^l$ or $P^r$ for some class $P$. 
Each of these classes  can be expressed as a Quillen negation of the class consisting of one or two maps of sets of size $\leqslant 2$, 
and 6 of them also lie in the orbit of $\{\emptyset\lra\{\bullet\}\}$. We now list these classes in the same notation we employ to denote maps of finite topological spaces.
The subscripts indicate where points map to; we also try to make this visually obvious by placing a point above its image.
\begin{figure}
	\hskip-12pt\noindent\includegraphicss[1.204949]{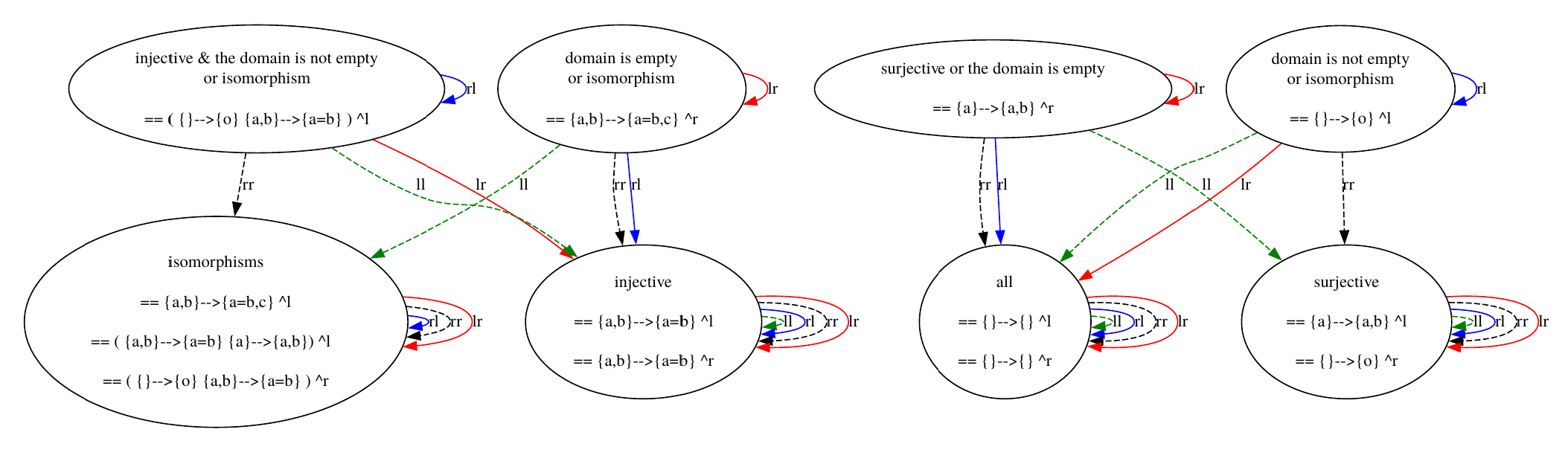}
	\caption{The Quillen negation monoid of the category of sets in the generators $ll,lr,rl,rr$. Quillen negations only. \label{Fig:QnSetsLLRR}} 
\end{figure}
\comment{\begin{figure}
	\hskip-12pt\noindent\includegraphicss[0.804949]{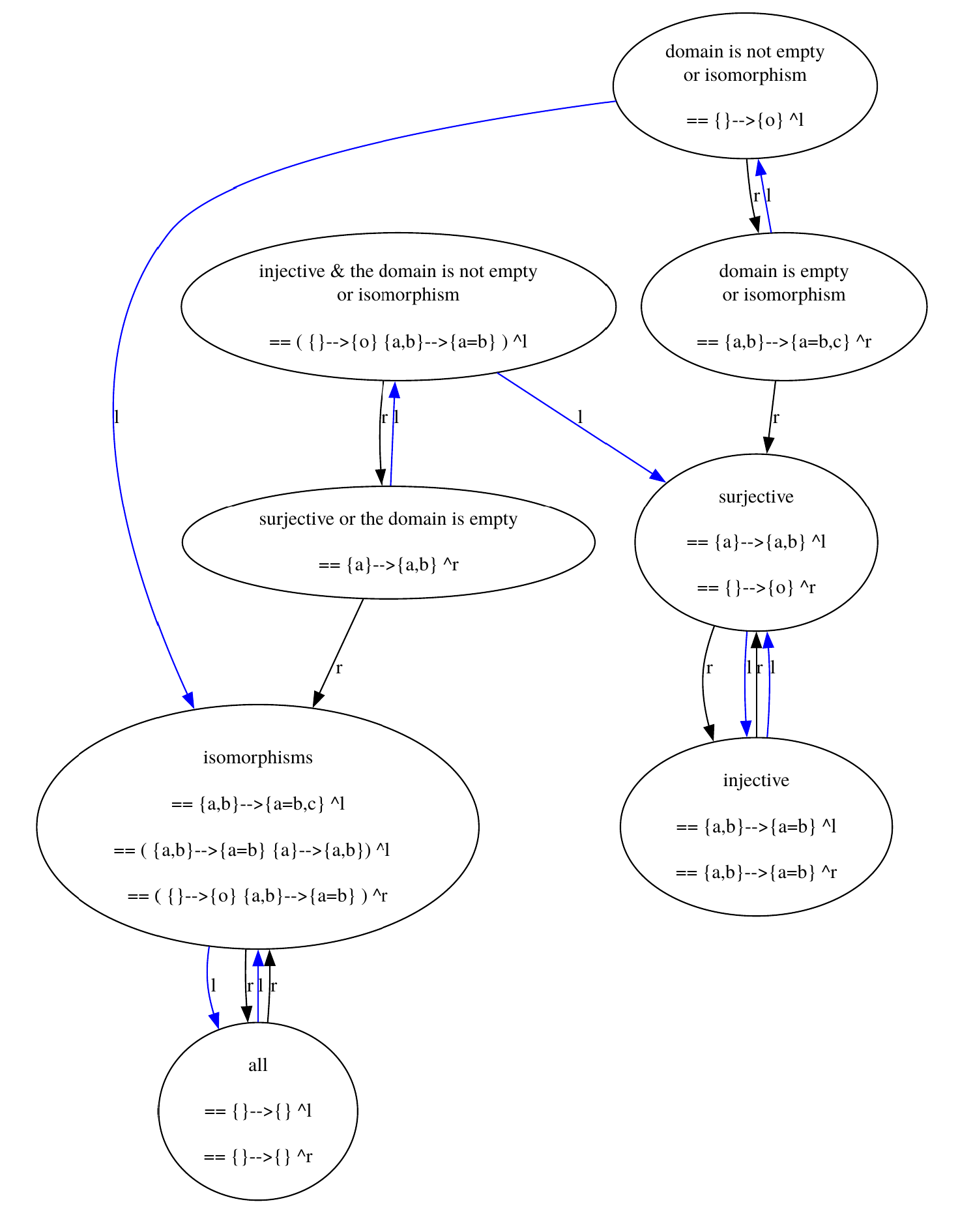}
	\caption{The Quillen negation monoid of the category of sets. Quillen negations only. \label{Fig:QnSetsNegationsOnly}} 
\end{figure}
}
$$\mathcalC^l=\text{(domain is not empty)}$$
$$\left({\displaystyle\bullet_a\phantom{\leftrightarrow\bullet_b}\atop{\displaystyle\Downarrow\atop \displaystyle \bullet_{a}\phantom{\leftrightarrow}  \bullet_b  }}\right)^l=
\mathcalC ^r=\mathcalC^{lrr}=\text{(surjections)}$$
$$\left({\displaystyle\bullet_a\,\, \bullet_b \atop{\displaystyle\Downarrow\atop \displaystyle \bullet_{a=b} }}\right)^l= 
\left({\displaystyle\bullet_a\,\,  \bullet_b \atop{\displaystyle\Downarrow\atop \displaystyle \bullet_{a=b} }}\right)^r=\mathcalC^{rl}=\text{(injections)}$$
$$\left({\displaystyle\bullet_a,  \bullet_b \phantom{,\bullet_c}  \atop{\displaystyle\Downarrow\atop \displaystyle \bullet_{a=b},\bullet_c }}\right)^l
=\left(\mathcalIso\right)^{lr}=\mathcalC^{ll}=\text{(isomorphisms)}$$
%
$$\left(\mathcalIso\right)^{l}=\left(\mathcalIso\right)^r
=\mathcalC^{lll}=\mathcalC^{llr}=\text{(all morphisms)}$$
$$\left({\displaystyle\bullet_a,  \bullet_b \phantom{,\bullet_c} \atop{\displaystyle\Downarrow\atop \displaystyle \bullet_{a=b},\bullet_c }}\right)^r
=\mathcalC^{lr}=\text{(the domain is empty or isomorphism)}$$
$$
\left( \mathcalCnobrackets, 
{\displaystyle\bullet_a\,\,\bullet_b \atop{\displaystyle\Downarrow\atop \displaystyle \bullet_{a=b} }}\right)^l=
 \text{(injective \& the domain is not empty or isomorphism)}$$
$$
\left(  
{\displaystyle\bullet_a	\phantom{\,\,\bullet_b} \atop{\displaystyle\Downarrow\atop \displaystyle \bullet_{a} \,\,\bullet_b }}\right)^r
  =\text{(surjective or the domain is empty)}$$

\section{Calculating the Quillen negations/orthogonals of maps of finite spaces
\label{sec:negations_of_maps}}
\subsection{Orthogonal classes, of classes consisting of a single morphism of finite spaces}
A few of the iterated negations of $\emptyset\lra\ptt$ are themselves Quillen negations of maps of finite spaces of small size, mostly left but sometimes right,
which define notions 
such as {\em surjective}, {\em injective}, {\em induced topology}, {\em Separation Axioms $T_0$ and $T_1$}, {\em subspaces}, and {\em closed subsets}. 

In this section we calculate those Quillen negations. 
The proofs amount to spelling out verbally in the usual language of open and closed subsets the meaning of the lifting property diagrams:
you treat the points in the finite spaces as {\em names} of subsets (namely, their preimages, in a left Quillen negation) or of points (namely,
their images, in a right Quillen negation). Once spelled out, the lifting property usually becomes equivalent to one of the standard basic definitions
in topology. 

Let us say the same in more detail. For left Quillen negations, we use that
to give a map to a finite topological space $Y$  
is the same as to  
give names to several open/closed subsets and specify their properties 
determined by the shape/combinatorics of $Y$. 
For right Quillen negations, we use that
to give a map from a finite topological space $A$  
is the same as 
to give names to several points and specify which are required to belong to the closure of which. 

\subsection{A summary of results: Definitions by example}
All of the Theorems in this section follow the same pattern: 
a property is defined with help of a simple example not having the property. 
Let us now briefly state our results in this terms. 
The table Fig.~\Comment{\ref{Fig:Table}}{\ref{Fig:Table2}} gives a summary of our results explained in this section.
and \S\ref{subsubsec:visuals} explains our notation for maps of finite spaces we use here,
although in most cases the reader might find it visually apparent. 

Surjectivity is both the left and the right Quillen negation 
of the two simplest non-surjections $\mathcalBnobrackets$ and $\mathcalCnobrackets$;
by this we mean that 
$\mathcalB^l=\mathcalC^r$  is the class of {\em surjections}. 
Similarly, injectivity is Quillen negation of the simple non-injections $\mathcalEnobrackets$ and $\mathcalAnobrackets$, 
i.e. $\mathcalE^l=\mathcalA^r$ is the class of {\em injections}. 
Examples $\{\star\leftrightarrow\star\}$ and $\Otoc$ of a non-$T_0$ and a non-$T_1$ space lead
to the definition of 
the classes of maps with $T_0$, resp.~$T_1$, fibres,
as $\mathcalE^r$, resp. $\mathcalDbrackets^r$.

The class $\mathcalDbrackets^l$ of maps where the topology on the domain is induced from the codomain, 
is defined with help of the simple map failing this property. 
An example of a non-closed non-injective map $\closedsubspace$ leads to the definition of
the class $\left(\closedsubspace\right)^l$ 
of {\em closed injections}, or, equivalently, {\em closed subspaces}. 

For the map $\mathcalDD$ the codomain fails to be the quotient of the domain,
and this leads to the class
$\mathcalDDbrackets^l$ of maps such that the codomain is  {\em the quotient of the domain  disjoint union with a discrete space}.

The map $\nondenseimage$ is the simplest example of a closed subset, and of a map failing to have dense image.
Accordingly, $\nondenseimagebrackets^l$ is the class of maps with dense image, and 
$$\nondenseimagebrackets^{lr}=\left(\closedsubspace\right)^l$$ is the class of closed subspaces.

The map $\mathcalGnobrackets$ is an example of a map $f:X\lra Y$ such that $\pi_0(f):\pi_0(X)\lra\pi_0(Y)$ is not surjective. 
Accordingly, whenever $\pi_0(Y)$ is finite, 
$\pi_0(f):\pi_0(X)\lra\pi_0(Y)$ is  surjective iff
$$ {{\displaystyle A \atop \displaystyle \phantom{f} \downarrow f }\atop {\displaystyle B}}\in \mathcalG^l$$

\subsection{Connected, injective, surjective, separation axiom $T_0$\label{subsec:conn}}

\subsubsection{$\protect\mathcalA^l$ is the class of maps such that $\pi_0(f):\pi_0(X)\lra\pi_0(Y)$ is injective, for $\pi_0(X)$ and $\pi_0(Y)$ finite.}
$\{\bullet,\bullet\}$ is perhaps the simplest example of a space failing to be connected.
\begin{theorem}[$\pi_0$-injective]\label{surjp0} Let $f:X\lra Y$ be such that $\pi_0(X)$ and $\pi_0(Y)$ are finite. Then 
	$f \in \mathcalA^{l}$ iff the induced map $\pi_0(f):\pi_0(X)\lra \pi_0(Y)$ is injective.
	In particular, $X$ is connected iff $ {\displaystyle X\atop{\displaystyle\downarrow\atop\displaystyle \ptt} } \in  \mathcalA^{l}$.
\end{theorem}
\begin{proof} 
	Each map from a space $X$ to a discrete space factors via $\pi_0(X)$, hence the commutative square in $X\lra Y \,\rtt\, \{\bullet,\bullet\}\lra\{\bullet\}$ fits into the diagram 
$$\xymatrix@C=2.39cm{X \ar[d] \ar[r] &\pi_0(X) \ar[d]\ar[r] & \{\bullet,\bullet\}=\pi_0( \{\bullet,\bullet\}) \ar[d] \\
	Y \ar[r] & \pi_0(Y) \ar[r] & \{\bullet\} }=\pi_0(\{\bullet\})$$
$\implies$:	If $\pi_0(f): \pi(X)\lra\pi_0(Y)$ is not injective, pick two $x_1\neq x_2\in \pi_0(X)$ with $\pi_0(f)(x_1)=\pi_0(f)(x_2)$.
	Then there is no lifting $Y\lra\{\bullet\}$ for any map $\pi_0(X)\lra\{\bullet,\bullet\}$ separating $x_1$ and $x_2$,
	and such a map exists because we assumed $\pi_0(X)$ to be finite and therefore discrete. 
$\impliedby$:
	If $\pi_0(f): \pi_0(X)\lra\pi_0(Y)$ is injective, 
	the lifting $\pi_0(Y)\lra \{\bullet,\bullet\}$ gives the required lifting $Y\lra\{\bullet,\bullet\}$.
\end{proof}
\begin{rema} A verification shows that $\Bbb Q_{\neq 0} \lra \Bbb Q $ does not have the left lifting property with respect to $\mathcalAnobrackets$, 
	yet the induced map from connected components of $\Bbb Q_{\neq 0}$ into that of $\Bbb Q$ is injective. This shows that 
	we do need an assumption on $\pi_0(X)$ or $\pi_0(Y)$. 
\end{rema}


\subsubsection{$\protect\mathcalA^r$ is the class of injections} Note that the map $\mathcalAnobrackets$ is an archetypal example of a non-injective map. 
The proof is a trivial calculation we spell out in detail to demonstrate our conventions on the notation of commutative diagrams.

\begin{theorem}\label{inj}\label{injR}
$\mathcalA^r$ is the class of injective maps.   
\end{theorem}
\begin{proof}
$\implies$:	Let $g:A \rightarrow B$ be an arbitrary map in $\mathcalA^r$. Consider the diagram \eqref{eq:injR}. 

Pick two points $a,b\in B$ such that $g(a)=g(b)$. As our convention  suggests,
	define $f$ by  $f(\bullet_a)=a$ and $f(\bullet_b)=b$.
Then the square is commutative, and we have a lifting $h$. Therefore $a=f(\bullet_a)=h(p(\bullet_a))=h(p(\bullet_b))=f(\bullet_b)=b$. 
	Since we chose $a$  and $b$ arbitrarily, the map $g$ is injective.

$\impliedby$: If $g:A\lra B$ is injective then then the commutativity of the diagram requires that $f(\bullet_a)=f(\bullet_b)$. 
	In this case define $h(\bullet_{a=b}):=f(\bullet_a)=f(\bullet_b)$.
	\begin{equation}
\label{eq:injR}
\begin{tikzcd}
\bigl\{\bullet_a,\bullet_b\bigr\} \arrow[r, "f"] \arrow[d, "p"]
& A \arrow[d, "i"] \\
	\{\bullet_{a=b}\}  \arrow[r, "g"] \arrow[ur, dashed, "h"]
& B
\end{tikzcd}\end{equation}
\end{proof}

\subsubsection{$\protect\mathcalB^l$ is the class of surjections.}  Note that the map $\mathcalBnobrackets$ is an example of a non-surjective map. We spell out the proof in detail to demonstrate our convertions regarding 
our use of commutative diagrams. Further below similar proofs shall often be replaced by a hint consisting of 
a commutative diagram.
\begin{theorem} \label{sur2}
The left orthogonal class  $\mathcalB^l$ is the class of surjections.
\end{theorem}
\begin{proof}
    Let $(p:X \rightarrow Y) \in \mathcalB^l$ be an arbitrary morphism.
Consider the commutative square 
	$$\xymatrix@C=2.39cm{X \ar[r] \ar[d]_{\mathcalB^l\ni p} & {\{\star^X\}} \ar[d]^i\\
	Y \ar[r]|g \ar[ru]|h & \{ \star^{\operatorname{Im}p}\leftrightarrow \star^{Y-\operatorname{Im}p}\} }$$
Recall that by our conventions, superscripts indicate the indented preimages, and in this case it means
	that we define $g$ by $g\inv(i(\star^{\Imm p})):=\Imm p$. 
Now, $p\in \mathcalB^l$ implies that the unique map 
$h:Y\lra\{\star^X\}$
makes the lower triangle commute, which means that 
$\operatorname{Im}g\subset \operatorname{Im}i$. 
	Hence $Y-\operatorname{Im}p=\emptyset$ and thus $p:X\lra Y$ is a surjection, as required. 
    \\
    Conversely: let  $p:X\lra Y$ be a  surjection. Then for each $y\in Y$
	there is $x\in X$ such that $y=p(x)$. Hence $g(y)=g(p(x))=i(\star^X)=\star^{\Imm p}$
	and the unique map $h:Y\lra\{\star^X\}$ makes the diagram commutative.
\end{proof}

\subsubsection{$\protect\mathcalB^r$ is related to Axiom $T_0$.} Although the map $\mathcalB$ is used in our main calculation, 
its right Quillen negation is not. We include it for completeness. 

Define an equivalence relation on a topological space as follows: two points $x,y\in X$ are equivalent 
iff 
the induced topology on $\{x,y\}\subset X$ is antidiscrete, i.e. there is a continuous map 
$g:\{\filledstar_x\leftrightarrow\filledstar_y\}\lra X$ such that $x=g(\filledstar_x)$ and $y=g(\filledstar_y)$. 
Call the equivalence classes of this relation {\em topologically indistinguishable components} of $X$. 
\begin{theorem} \label{br}
The right orthogonal class     $\mathcalB^r$ is the class of maps $p:X\lra Y$ such that 
for 
each topologically indistinguishable component of $X$, its image is 
a topologically indistinguishable component of $Y$.
\end{theorem}
\begin{proof}
Consider the diagram
$$\begin{tikzcd}
	\{\star_a\} \arrow[r] \arrow[d, "i"]
& X \arrow[d, "p"] \\
\bigl\{\star_a \leftrightarrow \star_b \bigr\} \arrow[r, "g"] \arrow[ur, dashed, "h"]
& Y
\end{tikzcd}$$
\end{proof}

\subsubsection{$\protect\mathcalC^l$ is the class with empty domain.} 
Note that $\mathcalCnobrackets$ is an archetypal example of a map with the empty domain.  

\begin{theorem}[non-empty domain]\label{cl}
The left orthogonal class $\mathcalC^{l}$ defines the class of  maps 
	$A \rightarrow B$ such that $A \neq \emptyset$ or $A=B=\emptyset$. 
\end{theorem}
\begin{proof} There is no map from a non-empty set to the empty set, hence the commutative square exists 
iff $A=\emptyset$. 
\end{proof}

\subsubsection{$\protect\mathcalC^r$ is the class of {\em surjections}.} Note that $\mathcalCnobrackets$ is an archetypal example of a non-surjective map. 
	\begin{theorem}[surjective]\label{cr}\label{sur}\label{surj1} 
$\mathcalC^r$ is the class of surjective maps. 
\end{theorem}
\begin{proof} A point $b\in B$ is the same as an arrow $\{\bullet\}\lra B$, and it has to lift. 
\begin{equation}
\label{eq:sur}
\begin{tikzcd}
\emptyset \arrow[r] \arrow[d]
& A \arrow[d, "p"] \\
\{ \bullet_{b_0}^{a_0} \} \arrow[r, "g"] \arrow[ur, dashed, "h"]
& B
\end{tikzcd}
\end{equation}
\end{proof}

\subsection{Induced topology, subspace, injective, separation axiom $T_0$ and $T_1$}	

\subsubsection{$\protect\mathcalDbrackets^l$ is the class of maps $A\lra B$ such that {\em the topology on $A$ is induced from $B$}.}
The map $\mathcalD$ is perhaps the simplest map such that the topology on the domain is not induced from the codomain. 
\begin{theorem}[induced]\label{indt}
    $\mathcalDbrackets^l$ is the class of maps $p:X\rightarrow Y$ such that 
	the topology on $X$ is induced from $Y$.
\end{theorem}
\begin{proof} An open subset $U\subset X$ is the same as an arrow $X\lra \{\Otoc\}$, and 
and the lifting arrow $Y\lra \{\Otoc\}$ is the same as a subset $V\subset Y$ such that $U=V\cap X$. 
	\begin{equation}
   \begin{tikzcd}
X \arrow[r, "f"] \arrow[d, "p"]
	   & \bigl\{\bullet^U_V \rightarrow {\resizebox{5pt}{!}{$\blacksquare$}^{X-U}_{Y-V} } \bigr\} \arrow[d] \\
Y \arrow[r, "g"] \arrow[ur, dashed, "h"]
	   & \{\bullet\}
\end{tikzcd} 
	\end{equation}
\end{proof}
	
\subsubsection{$\protect\mathcalDbrackets^r$ is the class of maps with $T_1$ fibres.} Note that $\Otoc$ is the simplest example of a space not satisfying Separation Axiom $T_1$, i.e. such that not each point is closed. 
\begin{theorem}[$T_1$]\label{T1f}\label{Dr}
		  $\mathcalDbrackets^r$ is the class of maps such that each fibre is a $T_1$ space. 
\end{theorem} 
\begin{proof} To give a pair of points $x\in \cl y$ in a space $X$ is the same as to give a map 
$\{\Otoc\} \lra X$. 
$$\begin{tikzcd}
	\bigl\{\overset{\bullet_{x_0}}{} \searrow \overset{}{\smallblackbox_{x_1}} \bigr\} \arrow[r,"f"] \arrow[d, "\pi"]
& X \arrow[d, "p"] \\
\bullet_{y_0} \arrow[r, "g"] \arrow[ur, dashed, "h"]
& Y
\end{tikzcd}$$
\end{proof}

\subsubsection{$\protect\mathcalE^l$ is the class of {\em injective} maps.} The map $\mathcalEnobrackets$ is not injective. 
\begin{theorem}[injective]\label{inj2}\label{injL}\label{El} 
    $\mathcalE^l$ is the class of injective maps.
\end{theorem}
\begin{proof} To give a pair of points $x,y\in X$ is the same as to give a map $\{\bullet,\bullet\}\lra X$. 
\end{proof}

\subsubsection{$\protect\mathcalDEbrackets^l$ is the class of maps of form $A\subset B$.}

\begin{theorem}[subspaces]\label{subspaces}\label{DEl} 
	$$\mathcalDEbrackets^l=\mathcalDEmorphismbrackets^l=\mathcalEDmorphismbrackets^l$$ 
	is the class of subspaces, i.e. the class of 
	injective maps $A\lra B$ such that the topology on $A$ is induced from $B$.
\end{theorem}
\begin{proof} 
	By Theorems \ref{indt}(induced) and \ref{inj2}(injective) 
	$$\mathcalDEbrackets^{l} = \mathcalDbrackets^l\bigcap \mathcalEbrackets^l=\text{(induced topology)} \bigcap  \text{(injection)} = \text{(subspaces)}.$$
Because orthogonals are closed under retracts,	the equalities   $$\mathcalDEbrackets^l=\mathcalDEmorphismbrackets^l=\mathcalEDmorphismbrackets$$ follow
from the fact that both morphisms $\mathcalD$ and $\mathcalEnobrackets$ are retracts of
	$\mathcalDEmorphism$ and $\mathcalEDmorphism$,
	and, conversely, both $\mathcalDEmorphism$ and $\mathcalEDmorphism$ are retracts of the product $\mathcalD\times\mathcalEnobrackets$.
\end{proof}

\subsubsection{$\protect\mathcalEbrackets^r$ is the class of maps with $T_0$ fibres.} Note that $\{\star\leftrightarrow\star\}$ is the simplest example of a space failing to have Separation Axiom $T_0$, i.e. having distinct topologically indistinguishable points. 
\begin{theorem}[$T_0$]\label{T0f}
	$\mathcalEbrackets^r$ is the class of maps with $T_0$ fibres.
\end{theorem}
\begin{proof}
	The proof is similar to the proof of Theorem~\ref{T1f}($T_1$). 
\begin{equation}\label{diagram:T0f}
\begin{tikzcd}
\bigl\{\star \leftrightarrow \star \bigr\} \arrow[r,"f"] \arrow[d]
& X \arrow[d, "p"] \\
\{\star\} \arrow[r, "g"] \arrow[ur, dashed, "h"]
& Y
\end{tikzcd}
\end{equation}
\end{proof}

\subsection{Discrete, quotient, and disjoint union}

\subsubsection{$\protect\mathcalFbracketsSection^l $ 
is the class of maps of form $A\subset \textrm{Im}\, A \sqcup D$ where $D$ is discrete and the topology on $A$ is induced from the codomain} 
\begin{theorem}[discrete]\label{disj}\label{Fl}
 An injective map $f:X\lra Y$ is in $\mathcalFbrackets^l$ iff it is form 
$f:X \rightarrow X \sqcup D$ where $D$ is discrete space.

	More generally, 
$\protect\mathcalFbrackets^l $ 
is the class of maps of form $A\subset \text{Im} A \sqcup D$ where $D$ is discrete and the topology on $A$ is induced from its image $\operatorname{Im} A$
in the codomain. 
\end{theorem}
\begin{proof}$\implies$:  
	Note that map of Theorem~\ref{indt}(induced) 
	is a retract of the morphism under consideration, therefore 
	$\mathcalFbrackets^l \subset\text{(induced topology)}$.

	Take an arbitrary map $i:A \rightarrow B$ in $\mathcalFbrackets^l \bigcap \text{(injections)}$. 
For a point $b\in B\setminus A$,  consider the following diagram where $b$ is the preimage of $\star_b$:
$$\xymatrix@C=42.239pt{
	A \ar[d]_i \ar[r] & {\!\!\left\{{\displaystyle \filledstar_b\leftrightarrow\filledstar_\downarrow\hphantom{\filledstar\leftrightarrow\filledstar}}\atop{\hphantom{\displaystyle \star\star_b\leftrightarrow\star}\displaystyle \filledstar\leftrightarrow\filledstar^A_{B\setminus b}}\right\}} \ar[d] \\
B \ar[r] \ar@{-->}[ru] & 
	     \left\{{\scriptsize {\star_b\leftrightarrow\star\leftrightarrow\star_{B\setminus b} } }\right\} }
$$
Then $b$ has to be the preimage of $\bullet_b$ and therefore is open.
$\impliedby$:  Note  that $A\rightarrow A\bigsqcup D$ is a cobase change of 
	$\emptyset \rightarrow D$ along $\emptyset\lra A$. 
	As right orthogonals are closed under base change, 
	it is enough for us to find the lifting only for the second map. We do so using that $D$ is discrete and the map $\mathcalF$ is surjective.
	The verification that the lifting property holds for surjective maps such that the topology is induced, is similar to the proof of  Theorem~\ref{indt}(induced).  
\end{proof}



\subsubsection{$\protect\mathcalDDbracketsSection$ is the class of quotients disjoint union with a discrete space.}

\begin{theorem}[quotient disjoint union with a discrete space]\label{quotient}\label{DDl}
 A surjective map $f:X\lra Y$ is in $\mathcalDDbrackets^l$ iff the topology on $Y$ is the quotient topology.

An injective map $f:X\lra Y$ is in $\mathcalDDbrackets^l$ iff it is form 
$f:X \rightarrow X \bigsqcup D$ where $D$ is a discrete space.

More generally, 
$\mathcalDDbrackets^l$ is the class of maps of form 
	form 
	$X \rightarrow X/\approx \bigsqcup D$ where $D$ is a discrete space and  
	$X/\approx $ is a quotient of $X$.
\end{theorem}
\begin{proof} 
	Ponder the following diagrams \eqref{eq:DDl}. 
	The first diagram says that a subset $V\subset Y$ (the preimage of one of the points in $\{\star \leftrightarrow \star\}$)
is open iff 
	its preimage is open, which  for a surjective map $f:X\lra Y$ is precisely the definition of a quotient topology on $Y$ induced by $f:X\lra Y$.

The second diagram implies that each point not in the image is open. For an injective map $f:X\lra Y$ this means that 
 that $Y=X\sqcup D$ for a discrete space $D$ and the map $f:X\lra Y$ is an inclusion $X\lra X\sqcup D$. 
A verification shows that $X\lra X\sqcup D $ does have the required lifting property. 
	\begin{equation}\label{eq:DDl}	\xymatrix{ X \ar[r] \ar[d]|f & \left\{\raisebox{2pt}{$\bullet^{f^{-1}(V)}_{V}$}\searrow\raisebox{-2pt}{$\smallblackbox^{X\setminus f^{-1}(V)}_{V}$}\right\} \ar[d] \\
	Y \ar[r] \ar@{-->}[ru] & \{\star_V\leftrightarrow\star_{Y\setminus V}\} } 
	\hphantom{ops ops ops} 
	\xymatrix{ X \ar[r] \ar[d]|f & \left\{\raisebox{2pt}{$\bullet^{\{y\}}_{\{y\}}$}\searrow\raisebox{-2pt}{$\smallblackbox^{X}_{Y\setminus\{y\}}$}\right\} \ar[d]^{\hphantom{.....}y\notin \Imm X }  \\
	Y \ar[r] \ar@{-->}[ru] & \{\star_{\{y\}} \leftrightarrow\star_{Y\setminus\{y\}} \} } 
	\end{equation}
\end{proof}

\subsection{
$\pi_0$-surjective, and surjective or empty}

\subsubsection{$\protect\left(\mathcalGnobracketsSection\right)^l$ is the class of maps $f$ such that $\pi_0(f):\pi_0(X)\lra \pi_0(Y)$ is surjective, for $\pi_0(Y)$ finite} 
\begin{theorem}[$\pi_0$-surjective.] \label{surpi0}
The left orthogonal class $\mathcalG^l$ is the class of 
maps $p:X\lra Y$ such that the image $\Imm p$ intersects each non-empty clopen subset of $Y$.

	In particular, if $\pi_0(Y)$ is finite, this means precisely that the induced map 
	$\pi_0(p):\pi_0(X)\lra\pi_0(Y)$ is surjective. 
\end{theorem}
\begin{proof} A non-empty clopen subset disjoint from the image of $p:X\lra Y$ gives rise
to the diagram below. 
	\begin{equation}\label{eq:surpi0}    \begin{tikzcd}
X \arrow[r,"f"] \arrow[d, "p"]
		& \left\{ \bullet^X \right\} \arrow[d, "i"] \\
Y \arrow[r, "g"] \arrow[ur, dashed, "h"]
		& \left\{ \bullet_{Y''}^{X}  \;\;\bullet_{Y'}^{\emptyset} \right\}
	\end{tikzcd}\end{equation}
	To see that the first claim implies the second one, it is enough to note that if $\pi_0(Y)$ is finite,
	then any point of $\pi_0(Y)$ not in the image corresponds to a clopen subset not intersection the image of $f:X\lra Y$. 
\end{proof}

\begin{rema} The following example shows that it is necessary to assume that $\pi_0(Y)$ is finite, or at least discrete.  
Consider the inclusion $i: \Bbb Q_{>0} \lra \Bbb Q_{\geqslant 0}$. of positive rationals into non-negative rationals. 
	Evidently $ \Bbb Q_{>0} \xra i  \Bbb Q_{\geqslant 0} \,\rtt\, \{\bullet\}\lra\{\bullet,\bullet\}$, yet 
	$0$ is a connected component of $\Bbb Q_{\geqslant 0}$ which is not in the image. 
\end{rema}

\subsubsection{${A\atop {\downarrow\atop B}}\in \protect\left(\mathcalGnobracketsSection\right)^r$ iff either $A=\emptyset$ or $A\lra B$ is surjective.}
Note that $\mathcalGnobracketsSection$ is perhaps the simplest example of a non-surjective map with a non-empty domain. 
\begin{theorem}[surjective or empty domain]\label{gr}
    $\mathcalG^r$ is the class of surjective maps, and maps of form $\emptyset \rightarrow Y$.
In other words, a map is in $\mathcalG^r$ iff it is surjective or its domain is empty. 
\end{theorem}
\begin{proof}
The proof is a trivial calculation, see diagram~\eqref{eq:gr}. 
\begin{equation}\label{eq:gr}
    \begin{tikzcd}
\{\bullet^{x_0} \}\arrow[r,"f"] \arrow[d, "i"]
& X \arrow[d, "p"] \\
\{\bullet^{x_0}_{y_0} \; \bullet^{x_1}_{y_1}\} \arrow[r, "g"] \arrow[ur, dashed, "h"]
& Y
\end{tikzcd}\end{equation}
\end{proof}

\subsection{Closed subsets, and dense image\label{subsec:closed_subsets}}

\subsubsection{$\left(\protect\closedsubspace\right)^l$ is the class of closed inclusions (i.e. closed subsets)}

\begin{theorem}[closed subspaces]\label{closedsubspacesL} 
	$\left(\closedsubspace\right)^l= \nondenseimagebrackets^{lr}$ 
	is the class of closed subspaces, i.e. of closed injective maps.
\end{theorem}
\begin{proof} 
By Theorem~\ref{lem:denseL} $\nondenseimagebrackets^l$ is the class of maps with dense image.
Consider the lifting diagrams. In the first diagram evidently 
$\pi(B)=\pi(\clB(A))\subset\cl(\alpha(A))\subset \cl_YX=X$, hence the lifting exist.
In the second diagram, existance of the lifting implies that $X$ is a closed subset of $Y$.
\def\clB{\operatorname{cl_B}}	
$$\begin{gathered}\xymatrix@C=2.39cm{ A \ar[r]|\alpha\ar[d]|{\text{(dense image)}} & X \ar[d]|{({X\subset Y}\text{ closed})}    
\\ B \ar[r]|\pi \ar@{-->}[ru]  & Y  }
\xymatrix@C=2.39cm{ 
X \ar[r]|{\id}\ar[d] & X  \ar[d]    \\
	\operatorname{cl}_Y X  \ar[r] \ar@{-->}[ru]   & Y }
\end{gathered}$$

It is only left to consider the first orthogonal class. 
Note that maps   $\mathcalD$ and $\mathcalEnobrackets$ are retracts of $\closedsubspace$, hence  
	therefore by Theorems \ref{indt}(induced) and \ref{inj2}(injective) 
	$$\left(\closedsubspace\right)^{l} \subset  \text{(induced topology)} \cap  \text{(injections)} = \text{(subspaces)}.$$

	Now let $A\subset B$ be a subspace in $\left(\closedsubspace\right)^l$.
	To see that $A=\operatorname{cl_B}A$ is closed, 
consider the diagram 
	$$\xymatrix{A \ar[d] \ar[r] & \left\{\raisebox{2pt}{$\star^\emptyset\leftrightarrow\star^\emptyset_{B\setminus A}\leftrightarrow\star^\emptyset$} \searrow \raisebox{-2pt}{$\smallblackbox^A_{\operatorname{cl}_BA}$}\right\}\ar[d] \\
	B \ar[r] \ar@{-->}[ru]  &   {\{\star=\star_{B\setminus A} \leftrightarrow \star=\bblackbox_A\}}}$$
	where  the preimage of $\smallblackbox$ is $A$ for both  horizontal arrows. Then the lifting has to send $\operatorname{cl_B}A$ to $\smallblackbox$, 
	hence $\operatorname{cl_B}A\cap (B\setminus A)=\emptyset$, i.e. $\operatorname{cl_B}A=A$ as required.

Now let $A\subset B$ be closed. The following diagram shows how to find the lifting. 
	$$\xymatrix{A \ar[d] \ar[r] & \left\{\raisebox{2pt}{$\star^{X}_X \leftrightarrow \star^{Y=A\setminus(Z\cup X)}_{B\setminus (Z\cup X)}\leftrightarrow\star^{Z'=A\cap Z\setminus Z''}_{Z\setminus Z''} $} \searrow \raisebox{-2pt}{$\smallblackbox^{Z''}_{\operatorname{cl}_B(Z'')=Z''}$}\right\}\ar[d] \\
	B \ar[r]|\xi \ar@{-->}[ru]  &   {\{ \star_{B\setminus Z} \leftrightarrow \star=\bblackbox_Z\}}}$$
\end{proof}

	We need the following lemma to prove Theorems~\ref{crllrr}(rllrr)  and \ref{crllrrr}(rllrrr).
The condition on the map in the Lemma implies that in each fibre there is a point contained in each non-empty open subset of the fibre. 

Recall that a point of a topological space is called {\em generic} iff it lies in each non-empty open subset, 
or, equivalently, its closure is the whole space.

\begin{theorem}\label{lem:closedsubspacesLR} 
	Let $f:X\lra Y$ be a map with 
		a section $s:Y\lra X$ 
		such that $Y$ $s(y)$ is generic in the fibre $f\inv(y)$ with induced topology for each $y\in Y$, i.e.
		$$s(y)\in \bigcap\limits_{\substack{U\cap f^{-1}(y)\neq\emptyset\\ U \text { is open}}} U.$$
Then 
	$f\in \left(\closedsubspace\right)^{lr}$, i.e.  
	for each space $B$ and its closed subset $A$ it holds $(A\subset_{\text{closed}} B) \,\rtt\, f$. 
\end{theorem}
\begin{proof}
Let $A\subset B$ be a closed subset of $A$. 
Define a lifting $h:B\lra X$ by $h_{|A}=\tau$ and $h_{|B\setminus A}=s\circ \beta _{|B\setminus A}$. 
	Then both upper and lower triangles commute by construction of $h$,
and we only need to show that $h$ is continuous. 

Let $U\subset X$ be open. 
By the assumption on $s$ we have that $s^{-1}(U)=\{y\in Y: f^{-1}(y)\cap U\neq \emptyset\}=f(U)$ is open,
	and hence  $h^{-1}(U)=\tau^{-1}(U)\cup \beta^{-1}(s^{-1}(U))\setminus A)=
	\beta^{-1}(s^{-1}(U))\setminus (A\setminus \tau^{-1}(U)) $ 
is an open set without a closed set,  hence open.
	This shows that the lifting is continuous. \end{proof}


 \def\mathcalDensebrackets{\left(\mathcalDense\right)}
\def\mathcalDensebrackets{\left(\mathcalDense\right)}
\def\nondenseimage{\mathcalDense} 
\def\nondenseimagebrackets{\mathcalDensebrackets}


\subsubsection{$\left(\protect\closedsubspace\right)^{ll}=\left(\protect\mathcalDenseSection\right)^l$ is the class of maps with dense image.}
Recall that by Theorem~\ref{closedsubspacesL} $\left(\closedsubspace\right)^{ll}$ is the class of closed subspaces
and therefore $\mathcalDensebrackets\in \left(\closedsubspace\right)^{ll}$.
\begin{theorem}[
		dense image]\label{lem:denseL}
	$\left(\closedsubspace\right)^{ll}=\left(\protect\mathcalDenseSection\right)^l$ 
	is the class of maps with dense image.
\end{theorem}
\begin{proof}
$\impliedby$: 
Consider the diagram where the square on the right is a pull-back.  
\def\clB{\operatorname{cl_B}}	$$\xymatrix@C=2.39cm{ A \ar[r]|\alpha\ar[d]|{\text{(dense image)}} & X \ar[d]|{({X\subset Y}\text{ closed})} \ar[r] & \{\smallblackbox\} \ar[d]  \\ 
	B \ar[r]|\pi \ar@{-->}[ru] \ar@{-->}[rru] & Y \ar[r]& \{\raisebox{2pt}{$\bullet_{Y\setminus X}$}\rotatebox{-12}{\ensuremath{\to}} \raisebox{-2pt}{$\smallblackbox_X$}\} 
}$$
	If the image $\Imm A$ of $A$ is dense in $B$, then $\pi(B)=\operatorname{cl_Y}(\alpha(A))\subset \operatorname{cl_Y}(X)=X$, and the lifting exists.

\noindent$\implies$:
	Consider the diagram 
$$\xymatrix@C=2.39cm{ 
A \ar[r]\ar[d] & \operatorname{cl_B} \Imm(A\lra B) \ar[d] \ar[r]  &  \{\smallblackbox\} \ar[d] \\
	B \ar[r] \ar@{-->}[ru] \ar@{-->}[rru] & B \ar[r]  & \{\raisebox{2pt}{$\bullet_{B\setminus \clB \Imm A}$}\rotatebox{-12}{\ensuremath{\to}} \raisebox{-2pt}{$\smallblackbox_{\clB \Imm A}$}\} }$$
	If the image $\Imm A$ of $A$ is dense in $B$, then $\pi(B)=\operatorname{cl_Y}(\alpha(A)\subset \operatorname{cl_Y}(X)=X$, and the lifting exists.
	The lifting exists iff $B\subset \operatorname{cl_B} \Imm(A\lra B) $, i.e.~in other words $\Imm(A\lra B)$ is dense in $B$.
\end{proof}

\section{Calculating all the Quillen negations of $\mathcalCnobrackets$\label{sec:lrls}} 

Below we calculate case by case the orbit of $\mathcalCnobrackets$ shown in Fig.~\ref{Fig:Orbit}.
The labels in Fig.~\ref{Fig:Orbit} sketch the statements of all the theorems in our calculation below, 
and many readers may find it more efficient to start by looking at Fig.~\ref{Fig:Orbit} 
and refer to our calculation only if necessary.

\subsection{Discrete, subspace, section}

Recall that by Theorem~\ref{sur} $\mathcalC^r$ is the class of surjective maps. 

\subsubsection{rl$\protect\mathcalC=\left(\protect\mathcalDDSection, \protect\mathcalEnobrackets\right)^l$:
disjoint union with a discrete set.}

\begin{theorem}[rl$\protect\mathcalC$: discrete]\label{crl}
	$$\mathcalC^{rl}=\left(\mathcalDDSection, \protect\mathcalEnobrackets\right)^l=\left({\mathcalF},\mathcalEnobrackets\right)^r=$$ 
	is the class of maps $i:A\rightarrow A \bigsqcup D$ where $D$ is discrete space.   
\end{theorem} 
\begin{proof} $\implies$:
Note that maps $\mathcalDD$, $\mathcalE$ and $\mathcalF$ are surjective, 
	 hence  Theorems \ref{DDl}(quotient), \ref{inj2}(injective) and \ref{disj}(discrete) imply that 
	$$\mathcalC^{rl} \subset \mathcalEbrackets^l\cap \mathcalFbrackets^l \bigcap \mathcalEbrackets^l = \{A \rightarrow A \sqcup D\,: \,D\text { discrete}\}$$

\noindent$\impliedby$: A verification shows that maps of this form lift with respect to surjections. 
\end{proof}
\begin{rema}
	$\mathcalC^{rl}$ is the the closure of $\mathcalC$ under cobase changes and coproducts. 
	Indeed, a map $A\lra A \sqcup D$ is a cobase change (push forward) of $\emptyset\lra D$ along $\emptyset\lra A$, 
	and $\emptyset \lra D$ is $\emptyset \lra \coprod\limits_{d\in D} \ptt$ if $D$ is discrete.

	Note that each right Quillen negation is closed under cobase changes and coproducts  by \cite[Lemma 3.4]{H}. 
\end{rema}

\subsubsection{rr=$\protect\mathcalDEbrackets^l$:  subspace.}
Recall that by Theorem~\ref{sur} $\mathcalC^r$ is the class of surjective maps, and that
by Lemma~\ref{subspaces} $\mathcalDEbrackets^l$ is the class of subspaces. 
\begin{theorem}[rr: subspaces]\label{crr}
$\mathcalC^{rr}=\mathcalDEbrackets^l$ is the class of subspaces.
\end{theorem}
\begin{proof} 
$\implies$: 
	Note that the map in Theorem~\ref{inj}(injective) is  surjective, i.e.~in $\mathcalC^r$, 
hence 
	$\mathcalC^{rr} \subset \mathcalA^r=\text{(injections)}$. 
	Following notation of \cite[I\S5]{Bourbaki}, each injective map $f:X\lra Y$ has a canonical decomposition 
	$X \xra g f(X) \xra \psi Y$ where $\psi$ is the canonical injection of the subspace $f(X)$ into $Y$,
	and $g$ is the bijection associated with $f$. Now consider the diagram 
$$\begin{tikzcd}
X \arrow[r, "Id"] \arrow[d,"g"]
& X \arrow[d, "f"] \\
f(X) \arrow[r, "\psi"] \arrow[ur, dashed, "Id"]
& Y
\end{tikzcd}$$
	The continuity of the unique lifting means that the topology on $X$ is induced and $g$ is an isomorphism. 
%
%

\noindent$\impliedby$: 
Let $i:X\subset Y$ be a subspace, and consider a diagram 
$$\begin{tikzcd}
A \arrow[r, "u"] \arrow[d, "p"]
& X \arrow[d, "i"] \\
B \arrow[r, "d"] \arrow[ur, dashed, "h"]
& Y
\end{tikzcd}$$
	If the square commutes, then $\Imm d \subset X$ and we may define the lifting by $h:=d_{|X}$.
\end{proof}

\subsubsection{lr=$\protect\mathcalC^{lr}$:  $\emptyset\lra B$.}
Recall that by Theorem~\ref{cl} $\mathcalC^l$ defines the class of maps $A\lra B$ where $A\neq \emptyset$. 
\begin{theorem}[lr: $\emptyset\lra B$]\label{clr}
$\mathcalC^{lr}$ is the class of maps of form $\emptyset \rightarrow B$.    
\end{theorem}
\begin{proof}   
$\implies$: For arbitrary $f:A \rightarrow B$ in $\mathcalC^{lr}$ such that 
	$A \neq \emptyset$ it holds that $f\in \mathcalC^l$, i.e.~$A\lra B\,\rtt\,\emptyset\lra Y$ for each $Y$. 
	Hence 
	$f\in \mathcalC^{lr}\cap \mathcalC^l=(isomorphisms)$. 

\noindent$\impliedby$: 
For $A=\emptyset$ there is no commutative square with sides $f$ and $\emptyset\lra\{\bullet\}$, hence the lifting property holds vacuously.
\end{proof}

\subsubsection{ll: isomorphisms.}
\begin{theorem}[ll: isomorphisms] \label{cll}
	$\mathcalC^{ll}$ is the class of isomorphisms. 
\end{theorem}
\begin{proof}
Note that maps  $\mathcalBnobrackets$, $\mathcalD$ and $\mathcalEnobrackets$ are in $\mathcalC^l$, 
	therefore by Theorems \ref{sur}(surjective), \ref{indt}(induced) and \ref{inj2}(injective)) 
$$\mathcalC^{ll} \subset \text{(surjections)} \cap \text{(induced topology)} \cap  \text{(surjection)} = \text{(Isomorphisms)}.$$
As any isomorphism is contained in any left or right orthogonal class,
and this implies the equality and completes the proof.
\end{proof}

\subsubsection{lrr is the class of maps having a section.}
Recall that $\mathcalC^{lr}$ defines the class of maps with empty domain. 
The following theorem is immediate.
\begin{theorem}[lrr: section]\label{clrr}
$\mathcalC^{lrr}$ is the class of maps having a section.    
\end{theorem}

\subsubsection{rrrl=rr} 
Recall that by Theorem~\ref{crr} $\mathcalC^{rr}$ is the class of subspaces. 

\begin{theorem}[rrrl=rr: subspaces]\label{crrrl} 
	$\mathcalC^{rrrl}=\mathcalC^{rr}$ is the class of subspaces.
\end{theorem}
\begin{proof} By Theorem~\ref{DEl} the class of subsets is a left Quillen negation $P^l$,  namely
	$$\mathcalC^{rr}=\mathcalDEbrackets^l  .$$
	Applying $^{rl}$ to both sides of the equation using the identity $lrl=l$  
	we get 
	 $$
	 \mathcalC^{rrrl}=	\mathcalDEbrackets^{lrl}=\mathcalDEbrackets^{l}=\mathcalC^{rr}$$
\end{proof}

\subsection{Complete lattice, separation axiom $T_1$}

\subsubsection{rrr is related to the class of complete lattices} 

Unfortunately, we were unable to give a nice description of the class $\mathcalC^{rrr}$. However, the following is sufficient to calculate its orthogonals.  
\begin{theorem}[rrr: lattice]\label{crrr} 
Each map in $\mathcalC^{rrr}$ is a quotient map admitting a section. Moreover, 
	each fibre is a retract of a Cartesian power of $\{\Otoc\}
	\times \{\star\leftrightarrow\star\}$. 

	A map $X\lra\ptt$ is in $\mathcalC^{rrr}$  iff $X$ is a retract of a Cartesian power of $\{\Otoc\}
	\times \{\star\leftrightarrow\star\}$.

	In particular, a partial  order $P$ is a complete lattice iff $P\lra\ptt$ is in  $\mathcalC^{rrr}$.
\end{theorem}
\begin{proof} By Lemma~\ref{lem:LatticesasRetracts} below a partial  order $P$ is a complete lattice iff it is a a retract of a Cartesian power of $\Otoc$ 
in the category of partial orders. 
	As the category of partial orders is a full subcategory of the category of topological spaces,  
	this shows that the last claim follows from the previous ones.

	Let $X \xrightarrow g Y $ denote an arbitrary map in $\mathcalC^{rrr}$. Each map $\emptyset \lra Y$ is in $\mathcalC^{rr}$, hence 
	$\emptyset\xra{\mathcalC^{rr}} Y\,\rtt\, X \xra g Y$  and thus $g$ has a section $s:Y\lra X$. A standard argument shows that each map admitting a section
	is necessarily a quotient map.
%
	By Theorem~\ref{DEl} $\mathcalC^{rr}=(F\lra \ptt)^l$ for a space $F=\{\star\leftrightarrow\star\rightarrow \smallblackbox\}$, and 	
	therefore by Lemma~\ref{XtoOretracts} $X\lra \ptt$ belongs to ${\mathcalC^{rrr}}=(F\lra\ptt)^{lr}$ 
	iff $X$ is a retract of a (possibly infinite) Cartesian power of $F$.

\end{proof}
%
The following fact is well-known, e.g. \cite[p.203]{RW} and \cite[p.295]{PR} mention this fact as standard. 
\begin{lemm}\label{lem:LatticesasRetracts} In the category of partial orders and monotone maps, the following are equivalent
	for a partial order $P$:
\begin{enumerate} \item Partial order $P$ is a complete lattice.
		\item Partial order $P$ is  a retract of a product of $\Otoc$.
		\item ${P\atop{\downarrow\atop\ptt}}\in \left(\mathcalD\right)^{lr}$
\end{enumerate}
\end{lemm}
\subsubsection{rrrr=$\protect\mathcalDbrackets^r$: the fibre is $T_1$.}

\begin{theorem}[rrrr=$\protect\mathcalDbrackets^r$: the fibre is $T_1$]\label{crrrr} 
$\mathcalC^{rrrr}=\mathcalDbrackets^r$ is the class of maps $f:X\lra Y$ such that  each fibre $f^{-1}(y)$, $y\in Y$, satisfies separation axiom $T_1$.
\end{theorem}
\begin{proof}$\implies$: By Theorem~\ref{crrr}  $\mathcalD\in\mathcalC^{rrr}$, and thus by Theorem~\ref{T1f}($T_1$) each fibre of any map in $\mathcalC^{rrrr}$ 
	satisfies separation Axiom $T_1$.
	
\noindent$\impliedby$: 
In the commutative square each fibre has to go to a single point because each fibre of $f$ is $T_1$,
hence the map lifts by the universal property of quotient maps. 
\end{proof}

\subsubsection{rrr.rr: isomorphisms.}
\begin{theorem}[rrr.rr: isomorphisms]\label{crrrrr} 
$\mathcalC^{rrr.rr}$ is the class of isomorphisms.
\end{theorem}
\begin{proof} By Theorem~\ref{crrrr}(rrr.r:$T_1$)  $\mathcalCnobrackets, \mathcalAnobrackets\in\mathcalC^{rrr.r}$, hence by Theorems~\ref{surj1}(surjections) and \ref{inj2}(injection) 
	$$\mathcalC^{rrrrr} \subset \mathcalC^r \cap \mathcalA^r = (\text{surjection})\cap(\text{injection})=(\text{bijections})$$
	and thus	each map in $\mathcalC^{rrr.rr}$ is a bijection. On the other hand, by Theorem~\ref{crrrr}($T_1$) 
	each bijection belongs to $\mathcalC^{rrrr}$, hence 
	each map in $\mathcalC^{rrr.rr}$ lifts against itself and thus is an isomorphism.
\end{proof}

\subsubsection{rrr.rl}
Unfortunately, we do not have a good description of $\mathcalC^{rrr.rl}$. 
The following suffices to calculate further Quillen negations. 

\def\quotientl{{\{\bullet\rightarrow\smallblackbox\}}\atop{\{\bullet\leftrightarrow\smallblackbox\}}}
\def\quotientl{{\bullet\rightarrow\smallblackbox}\atop{{\downarrow}\atop{\bullet \leftrightarrow \smallblackbox}}}

\begin{theorem}[rrr.rl]\label{crrrrl} 
	$\mathcalC^{rrrrl}$ is contained in the class of quotient maps.
	 $\mathcalC^{rrrrl}$  contains each quotient map such that no fibre has a non-trivial closed equivalence relation.                                                              
	
	In particular, the class $\mathcalC^{rrrrl}$ contains $\mathcalF$.
\end{theorem}
\begin{proof}$\implies$: By Theorem~\ref{DDl} and Theorem~\ref{El}  the class of quotient maps is a left orthogonal 
	$\left(\mathcalDD,\mathcalEnobrackets\right)^l$, and by Theorem~\ref{crrr}
	        $\mathcalC^{rrr}\subset (quotients)$. Hence 
$$\mathcalC^{rrrrl}\subset (\text{quotients})^{rl} = \left(\mathcalDD,\mathcalEnobrackets\right)^{lrl} = \left(\mathcalDD,\mathcalEnobrackets\right)^l = (\text{quotients})$$

\noindent$\impliedby$:	Now assume that $f:A\lra B$ 
	is a  quotient map such that no fibre has a non-trivial closed equivalence relation.                                                              
A space $F$ has a non-trivial closed equivalence relation iff there is a non-trivial map $F\lra F'$ to some space $F'$ satisfying Axiom $T_1$.
Hence 
the argument in the proof of Theorem~$\ref{crrrr}$ using the universal property of quotient maps applies, namely
	in each commutative square with sides $f$ and $g\in \mathcalC^{rrrr}$ each fibre of $f$ goes to a single point,
	and the lifting exists by the universal property of the quotient maps.
\end{proof}

\subsubsection{rrr.rll=rl} 

\begin{theorem}[rrr.rll=rl] 
	$\mathcalC^{rrr.rll}=\mathcalC^{rl}=\left(\mathcalF\right)^l$ is the class of maps of form $A\lra A\sqcup D$ where $D$ is a discrete space.
\end{theorem}
\begin{proof} By Theorem~\ref{Fl} $\left(\mathcalF\right)^l$ is the class of maps of form $A\lra A\sqcup D$ where $D$ is a discrete space.
	By Theorem~\ref{crrrrl} ${\mathcalF}\in \mathcalC^{rrr.rl}$, hence 
	$\mathcalC^{rrr.rll}\subset \left(\mathcalF\right)^l$.  
On the other hand, by Theorem~\ref{crl}  $A\lra A\sqcup D$ left-lifts with respect to any surjection, hence to each element of $\mathcalC^{rrr.rl}\subset(surjections)$. 
\end{proof}

\subsubsection{rrl=r: surjections.} The following theorem is immediate using that
that $P^{lrl}=P^l$ for each $P$, and that surjections are defined by a left Quillen negation. 

\begin{theorem}
$\mathcalC^{rrl}=\mathcalC^r$ is the class of surjection.    
\end{theorem}

\subsection{Connectivity, injectivity.} 

\subsubsection{rll=$\left(\protect\mathcalGnobracketsSection\right)^{l}$: $\pi_0(X)\lra \pi_0(Y)$ is surjective.} 
Recall that by 
Theorem~\ref{crl}(rl) 
	$\mathcalC^{rl}$ 
is the class of maps $i:A\rightarrow A \sqcup D$ where $D$ is discrete space.  
The following theorem is immediate from Theorem~\ref{surpi0} and the fact that each map in this class is 
a pull-back of the map $\mathcalGnobrackets$.
\begin{theorem}[rll: surjective $\pi_0$]  \label{CRLL}\label{crll}
$\mathcalC^{rll}=\left(\mathcalGnobrackets\right)^l$ is the class of maps such that the image intersects any non-empty clopen subset. 
	In particular, if $\pi_0(Y)$ is finite, this means that the induced map $\pi_0(X)\lra \pi_0(Y)$ is surjective.
\end{theorem}

\subsubsection{rlll=ll: isomorphisms}
\label{CRLLtest}
Recall that by Theorem~\ref{crll}[rll]  $\mathcalC^{rll}=\left(\mathcalGnobrackets\right)^l$.
\begin{theorem}[rlll: isomorphisms]\label{crlll}
$\mathcalC^{rlll}$ is the class of isomorphisms.
\end{theorem} 
\begin{proof} A verification 
	shows that the maps $\mathcalBnobrackets$,
	$\mathcalEnobrackets$, and $\mathcalD$ have the left lifting property with respect to $\mathcalGnobrackets$, i.e. belong to
	$\mathcalC^{rll}=\left(\mathcalGnobrackets\right)^l$, hence by
	Theorems~\ref{sur2}(surjective), \ref{injL}(injective), and
	\ref{indt}(induced) $$\mathcalC^{rlll}\subset \mathcalB^l \cap
	\mathcalEbrackets^l \cap \mathcalDbrackets^l=(\text{surjective})\cap
	(\text{injective}) \cap (\text{induced})$$ is the class of
	isomorphisms. 
\end{proof} 

\subsubsection{lrrr=$\protect\mathcalEbrackets^l=\protect\mathcalA^r$ is the class of injections}
	Recall that $\mathcalC^{lrr}$ is the class of maps admitting a section.
\begin{theorem}[lrrr: injective] \label{clrrr}
	$\mathcalC^{lrrr}=\mathcalA^r=\mathcalEbrackets^l$ is the class of injections. 
\end{theorem}
\begin{proof}$\implies$:
By Lemma~\ref{inj} and \ref{inj2} $\mathcalA^r=\mathcalEbrackets^l$ is the class of injections.
The map from Lemma \ref{inj}(injective) belongs to $\mathcalC^{lrr}$, hence 
	$\mathcalC^{lrrr}\subset \mathcalA^r=\text{(injections)}$. 

	$\impliedby$: Consider the commutative diagram \eqref{eq:15}  where $i:X\lra Y$ is injective and  $p:A\lra B$ in $\mathcalC^{lrr}$ has a section
	$s:B\lra A$ such that $p\circ s=\operatorname{id}_B$. Define the lifting by $h:=s\circ f$. Then the lower triangle commutes,
	and by injectivity this implies that the upper triangle commutes.
\begin{equation}\label{eq:15}
\begin{tikzcd}
A \arrow[r, "f"] \arrow[d, "p"]
& X \arrow[d, "i"] \\
B \arrow[r, "g"] \arrow[ur, dashed, "g"]
&  Y
\end{tikzcd}\end{equation}
\end{proof}

%
%
%
%

\subsubsection{rllr=$\left(\protect\mathcalGnobrackets\right)^{lr}$} Unfortunately, we do not have a complete description of this item.
The following is sufficient to calculate further orthogonals.

Recall that by Theorem~\ref{crll} $\mathcalC^{rll}=\left(\mathcalGnobrackets\right)^l$ 
	and that by Theorem~\ref{closedsubspacesL}(closed subspaces)
$\left(\closedsubspace\right)^{l}$ is the class of closed subspaces.
\begin{theorem}[rllr]\label{crllr} 
	$\mathcalC^{rllr}=\left(\mathcalGnobrackets\right)^{lr}$ contains 
 all closed subspace inclusions $X\lra Y$ such that \begin{itemize}\item 
			$X$ is the intersection of all clopen subsets containing it:
$X=\bigcap\limits_{X\subset U\subset_{\text{clopen}}Y} U$\end{itemize}
	and is contained in the class of closed inclusions (equivalently, closed subspaces).
	$$\mathcalC^{rllr}=\left(\mathcalGnobrackets\right)^{lr}\subset  \left(\closedsubspace\right)^l=(\text{closed subsets})$$

	In particular, it contains maps of form $X\longrightarrow X\bigsqcup Y$. 
\end{theorem}
\begin{proof} 
$\impliedby$: Consider the diagram 
	$$\xymatrix@C=2.29cm{ A \ar[r] \ar[d]_{\left(\mathcalGnobrackets\right)^l} &X \ar[d]|i \ar[r] & \{\bullet^X\} \ar[d]  \\ B \ar[r] \ar@{-->}[ru] \ar@{-->}[rur] & Y \ar[r] & \{\bullet^X_U,\bullet_{Y\setminus U}\}}$$
	Let $X\lra Y$ be a closed inclusion satisfying the condition of the theorem. Since $X\subset Y$ is a subspace, the lifting $B\lra X$ exists if $\Imm(B\lra Y)\subset X$. 
Any clopen subset $U\supset X$ defines a commutative subdiagram of solid arrows such the square on the right is a pullback. 
	By assumption $A \lra B$ is in $\mathcalG ^l$, hence there is a lifting $B\lra\{\bullet\}$ (i.e. the unique map $B\lra\ptt$ is a lifting).  
The commutativity of the lower triangle for the lifting $B\lra\{\bullet\}$ means precisely that $X\subset U$. 
	Hence, $\Imm(B\lra X)\subset U$. As $U$ was chosen to be an arbitrary clopen subset of $Y$ containing $X$, we see that 
	$\Imm(B\lra X)\subset \bigcap\limits_{X\subset U\subset_{\text{clopen}} Y}=A$, and thus the lifting $B\lra X$ exists. 
	Hence, the map $X\lra Y$ is  in $\mathcalC^{rllr}$ as required.

\noindent$\implies$: By Theorem~\ref{crll} $(\text{surjections})\subset\mathcalC^{rll}$ and therefore $$\mathcalC^{rllr}\subset (\text{surjections})^r = (\text{subspace})$$
	Now let $X\subset Y$ be a subspace.  and let $i:X\lra Y$ be the corresponding inclusion map in $\mathcalG^{lr}=\mathcalC^{rllr}$.
	Let $\operatorname{cl}_Y X \subset Y$ be the closure of $X$ in $Y$.	
	Evidently $X\lra \operatorname{cl}_Y X \,\rtt\, \{\bullet,\bullet\}\lra\ptt$, hence $X\lra \operatorname{cl}_Y X \,\rtt\, X\lra Y$ and therefore
	$X=\operatorname{cl}_Y X$ as required.
\end{proof}

\subsubsection{rllrr: has a section picking a generic point in each fibre\label{subsubsec:rllrr}}
The proof of this theorem uses a set-theoretic cofinality argument 
which likely fails in many subcategories of topological spaces considered in algebraic topology or homotopy theory.

Recall that a point of a topological space is called {\em generic} iff it lies in each open subset, or, equivalently, its closure is the whole space.

Recall that by Theorem~\ref{crll} $\mathcalC^{rllr}=\left(\mathcalGnobrackets\right)^{lr}$.

\begin{theorem}[rllrr]\label{crllrr} 
%
%
	A map $X\xrightarrow f Y\in \mathcalC^{rllrr}$ iff 
	there is a section $s:Y\lra X$ of $f:X\lra Y$ such that $s(y)$ is a generic point of the fibre $f\inv(y)$ with induced topology, i.e.
	$$s(y)\in \bigcap\limits_{\substack{U\cap f^{-1}(y)\neq\emptyset\\ U \text { is open}}} U.$$

	In particular,  
each map $\mathcalC^{rllrr}$ has a section, is a quotient map, and each fibre has a generic point, 
	and 
	 $${\mathcalD}\,, {\mathcalEnobrackets} ,{\closedsubspace}\in \mathcalC^{rllrr}\subset \mathcalC^{rrrrl}$$
\end{theorem}
\begin{proof} The last three conditions on the map are immediately implied by existence of such a section,
	and in turn imply the inclusion by Theorem~\ref{crrrrl}(rrrrl). 

\noindent$\impliedby$: 
	By Theorem~\ref{lem:closedsubspacesLR} existence of such a section implies 
	the right lifting property with respect to closed inclusions, i.e. $f\in (\text{closed subspaces})^{r}$. By Theorem~\ref{crllr} 
	$\mathcalC^{rllr}\subset (\text{closed subspaces})$ hence $f\in (\text{closed subspaces})^r \subset  \mathcalC^{rllrr}$ as required. 

\noindent$\implies$:
Let us first show the implication for maps of form $F\lra\ptt$.
The general case is the same but needs slightly more cumbersome notation. 

	Let $F$ be a space such that $F\lra\ptt$ is in $\mathcalC^{rllrr}$. 
Let $\kappa$ be a cardinal of large enough cofinality. 
	Equip $F\cup \kappa$ with a topology generated by (i) open subsets $U\sqcup \{ \beta\in\kappa: \alpha<\beta<\kappa\}$
	and (ii) $\{\beta\in \kappa: \beta<\alpha\}$, 
where $U$ ranges among open subsets of $F$, and $\alpha\in\kappa$. 
	By construction the subsets $F\sqcup \{ \beta\in\kappa: \alpha<\beta<\kappa\}$ are clopen, and thus 
	$F$ is the intersection of all clopen subsets of $F\cup \kappa$ containing $F$, and 
	thus by Theorem~\ref{crllr} $F\lra F\cup\kappa$ is in $\mathcalC^{rllr}$.
	Now consider the diagram 
	$$\xymatrix@C=2.39cm{ F \ar[r]|{\operatorname{id}} \ar[d]_{\mathcalC^{rllr}} & F \ar[d]^{\mathcalC^{rllrr}} \\
	F\cup\kappa \ar[r] \ar@{-->}[ru]|h& \ptt } $$
	Assume the lifting $h:F\cup\kappa\lra F$ exists. Then the preimage of a non-empty open subset of $F$ is an open subset of form (i) and thus 
	contains 
	a final segment $\{\beta\in\kappa:\beta>\alpha\}$ for some $\alpha<\kappa$. In particular, 
	each cofinal subset of $\kappa$ intersects the preimage of each non-empty open subset of $F$. 
	If we choose the cofinality of $\kappa$ to be large enough, e.g. $\operatorname{cof} \kappa > \operatorname{card} F$, 
	then there is a point $x_0\in F$ 
	such that its preimage $h^{-1}(x_0)\setminus F$ 
	is cofinal in $\kappa$. Therefore its preimage intersects the preimage of each non-empty open subset, 
	and therefore $x_0$ is contained in each non-empty open subset of $F$. Thus $s:\bullet\mapsto x_0$ is the required section, 
	and this completes the proof for $F\lra\ptt$. 

	Now we do the general case. Let $X\lra Y$ be in $\mathcalC^{rllrr}$.
	Let $\kappa$ be a cardinal of large enough cofinality.
	Equip $X\cup \kappa\times Y$ with a topology generated by open subsets $U\sqcup \{ \beta\in\kappa: \alpha<\beta<\kappa\}\times V$
	        and $\{\beta\in \kappa: \beta<\alpha\}\times Y$,
		where $U\subset X$, resp.~$V\subset Y$, ranges among open subsets of $X$, resp.~$Y$, such that $f(U)\subset V$, and $\alpha\in\kappa$.
		Note that the induced topology on each fibre of $f:X\lra Y$ coincides with the topology defined above. 
		Same as before, by construction $X$ is the intersection of all clopen subsets of $X\cup \kappa\times Y$ containing $X$, and
		        thus by Theorem~\ref{crllr} $X\lra X\cup\kappa\times Y$ is in $\mathcalC^{rllr}$.
			        Now consider the diagram 
					$$\xymatrix@C=2.39cm{ X \ar[r]|{\operatorname{id}} \ar[d]_{\mathcalC^{rllr}} & X \ar[d]^{\mathcalC^{rllrr}} \\
					        X\cup\kappa\times Y \ar[r] \ar@{-->}[ru]|h& Y } $$
						Assume the lifting $h:X\cup\kappa\times Y\lra X$ exists.

						Then the preimage of a non-empty open subset $U$  of $X$ contains
							a subset of form  $U\cup\{\beta\in\kappa:\beta>\alpha\}\times V$ for some $\alpha<\kappa$ 
							and $f(U)\subset V$ both open.

	If  $\operatorname{cof} \kappa>\operatorname{card} X^Y$ 
	then there is a cofinal subset  $S\subset  \kappa$ such that $h(\alpha,y)=h(\beta,y)$ for each $\alpha,\beta\in S$ and $y\in Y$.

		Because the induced topology on each fibre $F=f^{-1}(y)$ of $f:X\lra Y$ coincides with the topology defined on $F$ defined above, 
		the argument above implies 
		$h(\alpha,y)$ lies in each open subset of $X$ intersecting $f^{-1}(y)$, for each $y\in Y$ and $\alpha \in S$.

		Hence, for $\alpha\in S$ we have that  $s:=h(\alpha,-):Y\lra X$ is a section of $f:X\lra Y$ such that for each $y\in Y$
		$$s(y)\in \bigcap\limits_{\substack{U\cap f^{-1}(y)\neq\emptyset\\ U \text { is open}}} U.$$
\end{proof}
\def\card{\operatorname{card}} 
\subsubsection{rll.rrr=rrrr: each fibre is $T_1$}
Recall that by Theorem~\ref{crrrr}($T_1$) $\mathcalC^{rrrr}$ is the class of maps with $T_1$ fibres.
\begin{theorem}[rll.rrr: the fibre is $T_1$]\label{crllrrr}
	$\mathcalC^{rll.rrr}= \mathcalC^{rrr.r}$ is the class of maps such that each fibre satisfies separation axiom $T_1$. 
\end{theorem}
\begin{proof}$\implies$: 
	By Theorem~\ref{crllrr}[rllrr]  $\mathcalD\in\mathcalC^{rllrr}$, hence $\mathcalC^{rllrr.r}\subset \mathcalDbrackets^r=(\text{each fibre is }T_1)$. 

	$\impliedby$:  Theorem~\ref{crllrr}[rllrr]  implies that any map in $\mathcalC^{rllrr}$ is a quotient map such that 
	there is no non-trivial closed equivalence relation on each fibre, 
	hence by Theorem~\ref{crrrrl}(rrr.rl) $\mathcalC^{rllrr}\subset \mathcalC^{rrr.rl}$ and therefore
	applying the right Quillen negation to the inclusion $\mathcalC^{rllrr}\subset \mathcalC^{rrr.rl}$ gives %
        $$\mathcalC^{rllrrr}\supset \mathcalC^{rrr.rlr}=\mathcalC^{rrrr}=(\text{each fibre is }T_1).$$ %
\end{proof}

\subsubsection{rll.rrl=$\left(\protect\closedsubspace\right)^l$: is the class of closed inclusions.} 

\begin{theorem}[rll.rrl: closed subsets]\label{crllrrl}\label{CClosedSubsets}
	$\mathcalC^{rll.rrl}= \left(\closedsubspace\right)^l = \nondenseimagebrackets^{lr}$ 
	is the class of closed subspaces.
\end{theorem}   
\begin{proof} The second equality follows from Theorem~\ref{closedsubspacesL}. 

	The map $\closedsubspace$ satisfies the condition of Theorem~\ref{crllrr}(rllrr) and thus 
	${\closedsubspace}\,\in\,\mathcalC^{rllrr}$, hence by Theorem~\ref{closedsubspacesL}(closed subspaces)  
	$$\mathcalC^{rllrrl}\subset \left(\closedsubspace\right)^{l}= (\text{closed subspaces}).$$

	By Theorem~\ref{lem:closedsubspacesLR} closed subspaces have the left lifting property with respect to maps having a generic section,
	and thus 
	$$\mathcalC^{rllrr}\subset \left(\closedsubspace\right)^{lr}=(\text{closed subspaces})^{r}$$
	hence applying the left Quillen negation gives
	$$\mathcalC^{rllrrl}\supset \left(\closedsubspace\right)^{lrl}=\left(\closedsubspace\right)^{l}=(\text{closed subspaces}).$$
\end{proof}


\subsubsection{rll.rrll=$\left(\protect\closedsubspace\right)^{ll}=\protect\nondenseimageSectionbrackets^l$ is the call of maps having dense image.}
Recall that by Theorem~\ref{crllrrl} $\mathcalC^{rllrrl}$ is the class of closed subsets. 
\begin{theorem}[rll.rrll: dense image]\label{crllrrll}
	$\mathcalC^{rll.rrll}= \left(\closedsubspace\right)^{ll}=\nondenseimagebrackets^l$ 
	is the class of maps with dense image. 
\end{theorem}
	\begin{proof}  Theorem ~\ref{crllrrl}(rllrrl) implies the first equality. Theorem ~\ref{lem:denseL} implies the second equality, 
		as well as the fact that this is the class of maps with dense image.
	\end{proof}

\subsubsection{rll.rrlll=ll: isomorphisms.}

\begin{theorem}[rll.rrlll: isomorphisms]\label{crllrrlll}
	$\mathcalC^{rll.rrlll}= 
		(\text{dense image})^l$ 
		is the class of isomorphisms.
\end{theorem}
\begin{proof} 
	Note that maps  $\mathcalB$, $\mathcalD$ and $\mathcalE$ have dense image and thus are in $\mathcalC^{rllrrll}$ by Theorem~\ref{crllrrll}(rllrrll). 
By Theorems \ref{sur}(surjective), \ref{indt}(induced) and \ref{inj2}(injective) 
	$$\mathcalC^{rllrrlll}  \subset \text{(surjections)} \cap \text{(induced topology)} \cap  \text{(injections)} = \text{(Isomorphisms)}.$$
\end{proof}

\subsubsection{rllrrllr=rllrrl: closed subspace}\label{crllrrllr}
Recall that by Theorem~\ref{crllrrl} 
$\mathcalC^{rllrrl}$ is a right orthogonal, hence the required equality is implied by the identity $rlr=r$.

\subsection{Disjoint union, quotient, induced topology, the fibre is $T_0$} 

%
%
%

\subsubsection{lrrl: disjoint union}%
Recall that by Theorem~\ref{clrr} $\mathcalC^{lrr}$ is the class of maps having a section.

\begin{theorem}[lrrl: disjoint union]\label{clrrl}  
$\mathcalC^{lrrl}$ is the class of maps of form $A\longrightarrow A\bigsqcup B$. 
\end{theorem}
\begin{proof}$\impliedby$:  Consider  a commutative square 
$$\begin{tikzcd}
A \arrow[r, "f"] \arrow[d, "i"]
& X \arrow[d, "p"] \\
A \bigsqcup B \arrow[r, "g_1 \bigsqcup g_2"] \arrow[ur, dashed, "h:=g_1\bigsqcup  g_2\circ s"]
& Y \arrow[u, "s"]
\end{tikzcd}$$
	and define the lifting as $h:=f \sqcup g_2\circ s: A\bigsqcup B \lra X$. 
		It is continuous by definition of the disjoint union. 
		Alternatively, a more category-theoretic approach is as follows. Notice that $A\lra A\sqcup B$ 
		is a cobase change of $\emptyset\lra B$ along $\emptyset\lra A$, hence by \cite[Lemma 3.2]{H} 
		it belongs to $\mathcalC^{lrrl}$ whenever $\emptyset \lra A$ does.
		By Theorem~\ref{clr} $\emptyset\lra A$  belongs to $\mathcalC^{lr}\subset \left(\mathcalC^{lr}\right)^{rl}=\mathcalC^{lrrl}$.

	\noindent$\implies$: Let $p:X\lra Y$ be in $\mathcalC^{lrrl}$. 
	Consider the obvious map $p':X\sqcup Y\lra Y$. It has a section $Y\lra X\sqcup Y$ and therefore there is a lifting $h:Y\lra X\sqcup Y$ 
		in diagram \eqref{eq:17.1}. Then $Y=h^{-1}(X) \sqcup h^{-1}(Y)$ is the disjoint union of the preimages of clopen subsets $X$ and $Y$ of $X\sqcup Y$. 
		Commutativity of the upper triangle implies that the topology on $h^{-1}(X)\subset Y$ is induced from $X\subset X\sqcup Y$, 
		and commutativity of the lower triangle implies that the topology on $h^{-1}(Y)\subset Y$ is induced also. 
\begin{equation}\label{eq:17.1}\begin{tikzcd}
X \arrow[r, "\id_{X}"] \arrow[d, "p"]
& X \bigsqcup Y \arrow[d, "p'"] \\
Y \arrow[r, "\id_Y"] \arrow[ur, dashed, "h"]
& Y
\end{tikzcd}\end{equation}
\end{proof}

\begin{rema}
$\mathcalC^{lrrl}$ is the class of maps obtained from $\mathcalGnobrackets$ from by base change.
	Indeed, to give a map $B\lra \{\bullet,\bullet\}$ is the same as to give a disjoint union $A\sqcup B'$, and 
	the base change along this map is $A\lra A\sqcup B'$. It is also obtained from the class $\mathcalC^{lr}$ of maps with empty domain 
	by cobase change. The following diagram illustrates this.
	$$\begin{gathered}\xymatrix@C=2.39cm{ \emptyset \ar[r]\ar[d]  & A \ar[d] \\ B\ar[r] & A\sqcup B } \phantom{AAAA} 
	\xymatrix@C=2.39cm{ A  \ar[r]\ar[d] & \ptt \ar[d] \\ A\sqcup B \ar[r] & \{\bullet,\bullet\} }\end{gathered}
	$$
\end{rema}
\subsubsection{lrrrl=$\protect\mathcalA^{rl}=\left({\protect\mathcalBnobrackets}, {\protect\mathcalDDSection}\right)^l$: quotient}%
This calculation is a standard universal property of quotient map with respect to injective maps.

Recall that by Theorem~\ref{clrrr} $\mathcalC^{lrrr}$ is the class of injective maps.  
\begin{theorem}[lrrrl: quotient]\label{clrrrl}
	$\mathcalC^{lrrrl} =\left({\protect\mathcalBnobrackets}, \protect\mathcalDDSection\right)^l = \mathcalA^{rl} $ is the class of quotient maps, i.e. surjective map $p:A\lra B$ such that the topology on $B$ is pushed forward from $A$.
\end{theorem}
\begin{proof}$\implies$: 
	The maps $\mathcalBnobrackets$ and $\mathcalDD$ 
	are injective, hence by Theorem~\ref{sur2}(surjective) and Theorem~\ref{quotient}(quotient) 
	$\mathcalC^{lrrrl}\subset \mathcalB^l \cap \mathcalDDbrackets^l 
	=\text{(surjective)}\cap\text{(quotient)}$,
	as required. 

\noindent$\impliedby$: 
	Consider the diagram 
$$\begin{tikzcd}
A \arrow[r, "f"] \arrow[d, "p"]
 & X \arrow[d, "i"] \\
 B \arrow[r, "g"] \arrow[ur, dashed, "h"]
 &  Y 
 \end{tikzcd}$$
	A possibly non-continuous lifting $h:B\lra X$ exists because $p:A\lra B$ is surjective and $i:X\lra Y$ is injective.
	It 
	is continuous due to the universal property of quotient 
	topology: $h$ is continuous if and only if $p\circ h=f$ is continuous, which we do know. 
\end{proof}
\subsubsection{lrrrll=rl: disjoint union with a discrete space}
Recall that by Theorem~\ref{clrrrl}(quotients)   $\mathcalC^{lrrrl}$ is the class of quotient maps, and that by Theorem~\ref{crl}(rl) 
$\mathcalC^{rl}$
is the class of maps of form $A\lra A\sqcup D$ where $D$ is discrete.
\begin{theorem}[lrrrll=rl]\label{clrrrll}
	$\mathcalC^{lrrrll} =\mathcalC^{rl}$ is the class of maps of form $A\lra A\sqcup D$ where $D$ is discrete.
\end{theorem}
\begin{proof}$\impliedby$: A verification shows that 
	any such map lift with respect to any surjection. Alternatively, by Theorems~\ref{clrrrl}(quotient) and ~\ref{cr}(surjective) 
	$\mathcalC^{lrrrl} \subset\mathcalC^{r}$ and therefore $\mathcalC^{lrrrll} \supset \mathcalC^{rl}$. 
	\noindent$\implies$: 
	The maps $\mathcalA$ and  $\mathcalF$ are quotient maps and therefore belong to $\mathcalC^{lrrrl}$, hence by Theorem~\ref{inj}(injective) and \ref{Fl} 
	$$\mathcalC^{lrrrll} \subset \mathcalEbrackets^l \cap \mathcalFbrackets^l=(\text{injective}) \cap (\text{quotient disjoint union with discrete})$$
as required.
\end{proof}

\subsubsection{lrrrr=$\protect\mathcalEbrackets^{lr}=\left({\protect\mathcalBnobrackets }, {\protect\mathcalD}\right)^l$: induced topology}%
This calculation is a standard universal property of induced topology with respect to injective maps.
Recall that by Theorem~\ref{clrrr}(injective)  $\mathcalC^{lrrr}$ is the class of injective maps.  

\begin{theorem}[lrrrr: induced topology and surjective]\label{clrrrr}
	$\mathcalC^{lrrrr} =\left({\protect\mathcalBnobrackets}, {\protect\mathcalD}\right)^l = \mathcalEbrackets^{lr} $ is the class 
surjective map $p:A\lra B$ such that the topology on $A$ is induced from $B$.
\end{theorem}
\begin{proof}$\implies$: 
	The maps $\mathcalCnobrackets$ is injective and thus $\mathcalCnobrackets\in \mathcalC^{lrrr}$ by  Theorem~\ref{clrrr}(injective), 
	 hence by Theorem~\ref{sur}(surjective)  
	$\mathcalC^{lrrrr}\subset \mathcalC^r  
	=\text{(surjective)}$.
Now let $f:X\lra Y$ be a surjective map in $\mathcalC^{lrrrr}$. Let $X_Y$ be the set of points of $X$ equipped with the topology induced from $Y$,
	and let $X\lra X_Y$ be the obvious map. It is injective (in fact, bijective), hence $X \lra X_Y \,\rtt\, X\lra Y$ (see the diagram below).
	The lifting $X_Y \lra X$ is necessarily an identity on points, and its continuity means precisely that the topology on $X$ is induced from $Y$.
$$\begin{tikzcd}
X \arrow[r, "\id"] \arrow[d, ""]
 & X \arrow[d, "i"] \\
X_Y \arrow[r, ""] \arrow[ur, dashed, "h"]
 &  Y 
 \end{tikzcd}$$

\noindent$\impliedby$: 
	Consider the diagram 
\begin{tikzcd}
A \arrow[r, "f"] \arrow[d, "p"]
 & X \arrow[d, "i"] \\
 B \arrow[r, "g"] \arrow[ur, dashed, "h"]
 &  Y 
 \end{tikzcd}
	A possibly non-continuous lifting $h:B\lra X$ exists because $p:A\lra B$ is surjective and $i:X\lra Y$ is injective.
	It 
	is continuous due to the universal property of quotient 
	topology: $h$ is continuous if and only if $p\circ h=f$ is continuous, which we do know. 
\end{proof}

\subsubsection{lrrrrr=$\protect\mathcalEbrackets^r$: the fibre is $T_0$.}
	Recall that by Theorem~\ref{clrrrr} $\mathcalC^{lrrrr}$ is the class of surjective maps such that the topology on the codomain is induced. 

\begin{theorem}[lrrrrr=$\protect\mathcalEbrackets^r$: the fibre is $T_0$]\label{clrrrrr} 
$\mathcalC^{lrrrrr}=\mathcalEbrackets^r$ is the class of maps $f:X\lra Y$ such that  each fibre $f^{-1}(y)$, $y\in Y$, satisfies separation axiom $T_0$.
\end{theorem}
\begin{proof}$\implies$: By Theorem~\ref{clrrrr}  $\mathcalEnobrackets\in\mathcalC^{lrrrr}$, and thus by Theorem~\ref{T0f}($T_0$) 
	each fibre of any map in $\mathcalC^{lrrrrr}$ 
	satisfies separation Axiom $T_0$.
	
\noindent$\impliedby$: 
	Consider a diagram where $f\in \mathcalC^{lrrrr}$ and $g\in\mathcalC^{lrrrrr}$: 
	$$\xymatrix@C=2.39cm{ A \ar[r]|t \ar[d]|f & X \ar[d]|g \\ B \ar[r] \ar@{-->}[ru] & Y }$$ 
Each fibre of $A\xrightarrow f B $ is antidiscrete because the topology on $A$ is induced from $B$, hence 
	it maps into a single point of $Y$ by the assumption that each fibre of $X\xrightarrow g Y$ satisfies Axiom $T_0$. 
	This defines  a possibly not continuous  lifting $h:B\lra X$. 
	Howewer, it is continuous by the universal property of induced topology.
	The lower triangle is commutative because the map $A\xrightarrow f B$ is surjective.
\end{proof}

\subsubsection{lrrrr.rl=lrrrr} By Theorem~\ref{clrrrr}(induced topology) the class $\mathcalC^{lrrrr}$ 
	is also a left orthogonal, hence the identity $lrl=l$ applies.

\subsubsection{lrrrrr.r=ll: isomorphisms} Recall that by Theorem~\ref{clrrrrr}($T_0$)  the class $\mathcalC^{lrrrrr}$ 
is the class of maps with $T_0$ fibres.
\begin{theorem}\label{clrrrrr} 
	The class $\mathcalC^{lrrrrrr}$ is the class of isomorphisms.
\end{theorem}
\begin{proof} The maps $\mathcalC$ and $\mathcalA$ have $T_0$-fibres, hence 
	$$\mathcalC^{lrrrrrr}\subset \mathcalC^r \cap \mathcalA^r = (\text{surjective)} \cap (\text{injective})= (\text{bijections})$$
	Now, any bijection also $T_0$ fibres, hence any map in $\mathcalC^{lrrrrrr}$ lifts with respect to any bijection, 
	hence it lifts with respect to itself, and thus is an isomorphism. 
\end{proof}


\subsubsection{The computation is complete.} A verification shows that we have finished classification of all the iterated Quillen negations of $\mathcalC$
	in the category of (all) topological spaces.

													\newpage
\end{document}